\documentclass[12pt,a4paper]{article}
\usepackage{graphicx,tikz-cd} % Required for inserting images
\usepackage{amssymb}
\usepackage{amsmath,amsthm}
\usepackage{latexsym, amsfonts,hyperref}
\allowdisplaybreaks
\newcommand{\be}{\begin{equation}}
\newcommand{\ee}{\end{equation}}
\newcommand{\bea}{\begin{eqnarray}}
\newcommand{\eea}{\end{eqnarray}}
\newcommand{\bean}{\begin{eqnarray*}}
\newcommand{\eean}{\end{eqnarray*}}
\newcommand{\brray}{\begin{array}}
\newcommand{\erray}{\end{array}}
\newcommand{\ben}{\begin{equation}{nonumber}}
\newcommand{\een}{\end{equation}{nonumber}}

\newtheorem{dfn}{Definition}[section]
\newtheorem{thm}[dfn]{Theorem}
\newtheorem{lmma}[dfn]{Lemma}
\newtheorem{ppsn}[dfn]{Proposition}
\newtheorem{crlre}[dfn]{Corollary}
\newtheorem{xmpl}[dfn]{Example}
\newtheorem{rmrk}[dfn]{Remark}
\newcommand{\bdfn}{\begin{dfn}}
\newcommand{\bthm}{\begin{thm}}
\newcommand{\blmma}{\begin{lmma}}
\newcommand{\bppsn}{\begin{ppsn}}
\newcommand{\bcrlre}{\begin{crlre}}
\newcommand{\bxmpl}{\begin{xmpl}}
\newcommand{\brmrk}{\begin{rmrk}}
\newcommand{\edfn}{\end{dfn}}
\newcommand{\ethm}{\end{thm}}
\newcommand{\elmma}{\end{lmma}}
\newcommand{\eppsn}{\end{ppsn}}
\newcommand{\ecrlre}{\end{crlre}}
\newcommand{\exmpl}{\end{xmpl}}
\newcommand{\ermrk}{\end{rmrk}}
\newcommand{\IC}{\mathbb{C}}

\newcommand{\IN}{{I\! \! N}}
\newcommand{\IZ}{\mathbb{Z}}
\newcommand{\bbc}{\mathbb{C}}

\newcommand{\id}{\text{id}}
\newcommand{\pol}{\mathrm{Pol}}
\newcommand{\q}{\mathcal{Q}}

\newcommand{\clb}{{\cal B}}
\newcommand{\clc}{{\cal C}}
\newcommand{\cld}{{\cal D}}

\newcommand{\clg}{{\cal G}}
\newcommand{\clh}{{\cal H}}

\newcommand{\clk}{{\cal K}}

\newcommand{\clq}{{\cal Q}}

\newcommand{\clu}{{\cal U}}

\def\a*{{\cal A}_{h,*}}
\def\B{{\cal B}(h)}
\def\B1{{\cal B}_1(h)}
\def\b{{\cal B}^{\rm s.a.}(h)}
\def\b1{{\cal B}^{\rm s.a.}_1(h)}

\newcommand{\ot}{\otimes}

\newcommand{\raro}{\rightarrow}

\def \qed {$\Box$}
\begin{document}
\[
\]
\begin{center}
{\large {\bf  Projective corepresentations and cohomology of compact quantum groups}}\\
by\\
{\large Prabuddha Chakraborty{\footnote{ Systems science and informatics unit, Indian Statistical Institute, Bangalore}}}\\
%{ Stat-Math Unit, Kolkata,}
{\large Debashish Goswami{ \footnote{Stat-Math Unit, Indian Statistical Institute, Kolkata}}
,
{\large Kiran Maity}\footnote{Stat-Math Unit, Indian Statistical Institute, Kolkata}},
%\footnote{ Support by a research fellowship from Indian Statistical Institute is gratefully acknowledged.}}},
\\
%{ Indian Statistical Institute}\\
%{ 203, B. T. Road, Kolkata 700 108, India}\\
{e-mail: prabuddha@isibang.ac.in, debashish\_goswami@yahoo.co.in, maitykiran01@gmail.com,}

\end{center}

\begin{abstract}
We study projective unitary (co)representations of compact quantum groups and the associated second cohomology theory. 
 We introduce left/right/bi/strongly projective corepresentations and study them in details. In particular, we prove that given any compact quantum group $\q$, there are   compact quantum groups $\tilde{\q_l}, \tilde{\q_r}, {\tilde \q}_{bi}, {\tilde \q}_{stp}$, each of which contains $\q$ as a Woronowicz subalgebra and every left/right/bi/strongly projective unitary corepresentation of $\q$ lifts to a linear corepresentation of  these quantum groups respectively. We observe that the strongly projective corepresentations are associated with the second invariant ($S^1$-valued) cohomology $H^2_{uinv}(\cdot)$ of the quantum group. We define a suitable analogue of normalizer of a compact quantum group in a bigger compact quantum group and using this, associate a canonical discrete group $\Gamma_\q$ to a compact quantum group $\q$ which is an alternative generalization of the second group cohomology and we show by an example that $\Gamma_\q$ in general may be different from $H^2_{uinv}(\q,S^1) $. We also illustrate with examples that $\Gamma_{\q}$ is more suitable than $H^2_{uinv}(., S^1)$ for labeling topological phases under quantum group symmetry.

\end{abstract}

\section{Introduction}
Symmetry plays a major role in both mathematics and physics. This triggered a flurry of research in the theory of groups and their representations, both in analytic and algebraic frameworks. In the formalism of quantum mechanics in terms of operators on Hilbert spaces, it is often natural to consider symmetry as a map on the level of rays, that is on the set of unit vectors up to scalar multiplication. This leads to consideration of projective representation of the symmetry group. Mathematically, the theory of projective representation of group is closely related with group extension and cohomology theory\cite{groupch}.

Quantum groups and more generally the theory of rigid tensor categories have generalized the classical group symmetry in mathematical physiscs. Beginning with the algebraic formulation due to Drinfield and Jimbo\cite{drin1},\cite{drin2},\cite{jim} which was motivated by questions in physics related to the solution of quantum Yang-Baxter equations, the theory of quantum groups and Hopf algebras have traversed a long way, with various analytic formalism due to a number of mathematician, most notably Woronowicz \cite{woronowicz1987compact}, Vaes-Kusterman \cite{kus1},\cite{kus2},\cite{kus3} Van Daele \cite{mulhopf},\cite{van2},\cite{van3},\cite{van4},\cite{van5}. In physics, there is a lot of interaction between the mathematical theory of quantum groups and tensor categories with the emerging field of topological states of matter and quantum computation\cite{chen,cen}.

This makes it a natural question whether one can extend the classical theory of projective representation of group to the realm of quantum groups. Substantial work in this direction have already been done by  Kenny De Commer(\cite{proj2DeMaNe}, \cite{proj1kenny}, \cite{galois}, \cite{galois1}, \cite{tannaka}, \cite{tannaka2}), Sergey Neshveyev(\cite{neshveyev2011second}), \cite{neshveyev2012autoequivalences}), \cite{ne2016classification}), Lars Tuset(\cite{neshveyev2011symmetric}), Makoto Yamashita(\cite{ne2016classification}). The goal of present article is to contribute to understanding of projective (co)representation of compact quantum groups and study the associated cohomology. Some of the main results obtained by us concern extension of projective corepresentations of a given (compact) quantum group to (linear) corepresentations of a bigger quantum group . Indeed, using Tannaka-Krein reconstruction theorem, for a given compact quantum group, we prove existence of a possibly larger compact quantum group such that any unitary projective corepresentation of the original quantum group lifts to a unitary corepresentation of the bigger quantum group. In fact, we carry out such enveloping construction for various types (left, right,bi) of projective corepresentation defined by us. For strongly projective corepresentations defined by us, which in a sense are the closet to the classical case, we can relaize the dual of the original quantum group as a normal quantum subgroup of (dual of) the corresponding envelope. This leads to connection with the the second invariant cohomology of quantum groups in the sense of (\cite{snbook}).We define a suitable analogue of normalizer of a compact quantum group in a bigger compact quantum group and using this, associate a canonical discrete group $\Gamma_\q$ to a compact quantum group $\q$ which is an alternative generalization of the second group cohomology and we show by an example that $\Gamma_Q$ in general may be different from $H^2_{uinv}(\q,S^1) $.  '

%In the last part of the thesis we compute such cohomology for a few quantum groups explicitly using the techniques of tensor category and fiber functors.

Let us briefly discuss some possible applications of our results to the domain of physics, more precisely topological phases. Some quantum systems permit a gapped spectrum with a single or degenerate vacuum state. In order to specify conditions for the latter, a natural definition of symmetry involves of maps that preserve the transition probability of rays. Changing the point of view from rays to specific states, the symmetry action becomes projective. Specifically, symmetries are operators on the Hilbert space that are either linear and unitary or anti-linear and anti-unitary. The first case can be interpreted as a projective representation of groups. The vacuum state degeneracy is then contigent to the projective phase, whenever the symmetry is preserved by the Hamiltonian in question. The issue, whether a vacuum is degenerate, is then solved by answering whether a phase redefinition on the states can eliminate the projective phase. This problem is related to the field of group cohomology. Different vacuum state are parametrized by different equivalence classes in the second group cohomology, $H^{2}(G,U(1))$ of the symmetry group $G.$ Symmetry protected topological phase are found to be exactly those vacuum states with non-trivial elements in $H^{2}(G,U(1)).$ When one replaces group symmetry on quantum mechanical system by a natural quantum group symmetry, it is clear that projective corepresentation and cohomology of quantum groups will play a crucial role to label and understand topologiacl phase. In this context, we will illustrate with some physical examples that the group $\Gamma_{\q}$ is in fact better suited to label topological phases arising from quantum group symmetry than $H^2_{uinv}(\q, S^1).$

\section{Preliminaries}
\subsection{Compact and discrete quantum groups}
\bdfn
 \label{fonts} 
 Let $\q$ be a unital $C^{\ast}$-algebra and $\Delta$ be a unital $C^{\ast}$-homomorphism from $\q$ to $\q \otimes \q.$ Then $(\q,\Delta)$ is said to be a compact quantum group if it satisfying the following properties:
 \begin{enumerate}
     \item [1)] $(\Delta\otimes \id)\Delta =(\id \otimes \Delta)\Delta,$
     \item[2)]  Each of the linear subspaces  $\Delta(\q)(1\otimes\q)$ and $\Delta(\q))(\q \otimes 1)$  is norm dense in $\q \otimes\q .$
 \end{enumerate}
     
\edfn
\bdfn\label{Haar}
    The Haar state $h$ on a CQG $\q$ is the unique state on ($\q$,$\Delta$) which satisfies the following conditions:
        $(\id \ot h)(\Delta (a))= (h \ot \id) \Delta(a)= h(a) 1_{\q}.$
\edfn

   The Haar state on an arbitrary CQG $\q$  is a non-commutative analogue of Haar measure on a compact group $G$.
   \bdfn 
A  Woronowicz $C^*$-subalgebra $\q_2$ of a compact quantum group $ ( \q_1, \Delta ) $ is a $ C^* $ subalgebra of $ \q_1 $ such that $\Delta(Q_{2})\subseteq Q_{2} \ot Q_{2}$ and the inclusion map $i:$ $ \q_2 \rightarrow \q_1 $ is a CQG morphism. 
 \edfn
 \bdfn  
A discrete quantum group (DQG) is a pair $(A,\Delta)$, where $A$ is $C_{0}$ direct sum of full matrix algebra  say $\overline{\bigoplus_{i\in I,\,C_{o}} M_{n_{i}}(\bbc)}$, where $I$ is some index set and $\Delta$ is a comultiplication such that $(A_{0}, \Delta/_{A_0})$ becomes a multiplier Hopf  $*$- algebra, where $A_{0}=\oplus_{i\in I} M_{n_i}(\bbc)$ (algebraic finite sum) without $C^{\ast}$ closure.
\edfn 

\brmrk If $A$ is a DQG then there exists a unique counit $\epsilon$ which is a $*$ homomorphism from $A$ to $\bbc.$ It will satisfy\begin{align*}
    (\epsilon\otimes Id)\Delta(a)=(Id\otimes \epsilon)\Delta(a)=a,
\end{align*} for all $a\in A.$
\ermrk

%{\bf Dual DQG of a CQG}

\vspace{1mm}

There exists a canonical duality between CQG and DQG in a way such that for a CQG $\q$ with complete enumeration of mutually inequivalent corepresentations $\{U_{i}:i\in I\}$, where $I$ is an indexed set and $U_{i}\in B(\clh_{i})\otimes\q$. Then the corresponding dual DQG is given by $\hat{\q}=\overline{\oplus_{i\in I} B(\clh_{i})}$, where we take the $C_{0}$ closure.

%\brmrk \begin{enumerate}
   % \item [1)] If $A$ is a finite dimensional DQG, then $A$ is also a compact quantum group.
  %  \item[2)] If $A$ is a finite dimensional DQG and $A=\oplus M_{n_{i}}(\bbc)$, then there is a $n_{i_{0}}$ such that $M_{n_{i_{0}}}=\bbc $ because $\epsilon $ is a non-zero homomorphism on each $M_{n_{i}}$ and we know $M_{n_{i}}$ is a simple algebra for $n_{i}>1$, so $\epsilon (M_{n_{i}})=0$ for  $n_{i}>1.$ But $\epsilon$ is a non-zero homomorphism, hence there is some  $n_{i_{0}}$ such that $n_{i_{0}}=1.$
%\end{enumerate}
%\ermrk
\bdfn\label{act}\cite{pod}

    Let $A$ be a $C^{\ast}$-algebra. A right coaction of a CQG $\q$ on $A$ is a  *-homomorphism
       $\delta$:$A$ $\to $$A$ $\otimes$ $\q$, which satisfies the following properties:
       
       \begin{enumerate}
           \item [(i)]$\delta$ intertwines the co-multiplication, meaning that 
       ($\delta$ $\otimes$ Id)$\delta$=(Id$\otimes$$\Delta$)$\delta$.
       \item[(ii)]$\delta$ satisfy the density conditions [$\delta$($A$)(1$\otimes$ $\q$)]=$A$$\otimes \q$. \end{enumerate}
 \edfn

      Similarly, we can define the left coaction of a CQG on a $C^{*}$ algebra $A.$
       
            \brmrk 
            %It is well known that condition (ii) is equivalent to the existence of a norm-dense, unital $\ast$-subalgebra $A_{0}$ of $A$ such that $\delta(A_{0}) \subseteq A_{0} \ot_{\rm alg} \pol(\q)$ and on $A_{0}$, $({\rm id } \otimes \epsilon) \circ \delta={\rm id}$. From this, we can conclude that
            If $A$ is a finite dimensional algebra then $\delta(A)\subseteq A\ot \pol(\q).$
            \ermrk
 \bdfn
            
        A right coaction $\delta$ of $\q$ on $A$ is said to be ergodic if {$A^\delta$}=\{$\delta($a$)$=$a$$\otimes$ 1\} is isomorphic to the $\bbc$, that means the fixed point set of $\delta$ is a one-dimensional algebra.
        \edfn

    % \bppsn 
       % Let $\q$ be a CQG. Let ($A$,$\delta$) be a unital $\q$-$C^*$ algebra. If $\delta$ is ergodic then there exist a unitary corepresentation $V_{A}$$\in$$M$($K$($L^2$($A$)$\otimes$ $\q$) of $\q$ such that $\delta($a$)$=$V_{A}$($a$$\otimes$ 1)$V_{A}^*$.
       % \eppsn 
     
       % A similar result can obtained when $\delta$ is not a Ergodic coaction. For that we have to use the conditional expectation $E$:$A$$\rightarrow$ $A^\delta$ given by $E(a)$=(id$\otimes$ $h$)$\delta($a$)$ for all $a$ in $A$, to make a Hilbert $A^\delta$ module $L^2$($A$,$E_\delta$).
        \subsection{Unitary tensor categories and Tannaka-Krein reconstruction theorem}
      \bdfn A category $\clc$ is called a $C^{\ast}$ category if\begin{enumerate}
    \item [i)] For all objects $U$ and $V, Mor(U,V)$ is a Banach space and the map \begin{align*}
        Mor(V,W)\times Mor(U,V)\to Mor(U,W), ~given~by~(S,T)\to ST,
    \end{align*} is a  bilinear map and the following holds: $\|ST\|\leq \|S\| \|T\|;$
    \item[ii)]There exists an antilinear contravariant functor $\ast:\clc\to\clc$ which maps object to the same object  and for a $T\in Mor(U,V),T^{\ast}\in Mor(V,U)$ and the followings hod :\begin{enumerate}
        \item [a)]$T^{\ast\ast}=T$, for $T\in Mor(U,V)$;
        \item[b)]$\|T^{\ast}T\|=\|T\|^{2}$. $End(U)=~Mor(U,U)$ is a unital $C^{\ast}$-algebra.
        \item[c)] For any $T\in Mor(U,V),$ the element $T^{\ast}T\in End(U)$  is a positive element.
    \end{enumerate}
\end{enumerate}\edfn

     \bdfn
     A $C^*$-category $\mathcal{C}$ is called a $C^*$-tensor category if it is a tensor category such that:

\begin{enumerate}
    \item The associativity and unit morphisms are unitary.
    \item The involution is compatible with the tensor product:
    \[
    (S \otimes T)^* = S^* \otimes T^*
    \]
    for all morphisms $S,T$.
    \item $\clc$ is closed under finite direct sums.  If we choose any two objects $U,V$, then there exists an object $W$ and isometries $u\in Mor(U,W)~and~v\in Mor(U,W)$ such that $uu^{\ast}+vv^{\ast}=1$;
\item If $U\in Obj({\clc})$ and $P\in End(U)$ is a projection then there exists an object $V$ and isometry $v\in Mor(U,V)$ such that $vv^{\ast}=P.$ The zero object is the corresponding subobject for the zero projection in $End(U).$
\item for the unit object $End(1)=\bbc1.$

\end{enumerate}
\edfn

\bdfn
Let $U \in \mathrm{Obj}(\mathcal{C})$. A conjugate of $U$ is an object $\overline{U}$ together with morphisms
\[
R : \mathbf{1} \to \overline{U} \otimes U, 
\quad 
\overline{R} : \mathbf{1} \to U \otimes \overline{U}
\]
satisfying
\[
(\overline{R}^* \otimes \mathrm{id})(\mathrm{id} \otimes R) = \mathrm{id}, 
\quad
(R^* \otimes \mathrm{id})(\mathrm{id} \otimes \overline{R}) = \mathrm{id}.
\]

A $C^*$-tensor category is \emph{rigid} if every object has a conjugate. Conjugates are unique up to isomorphism. Without loss of generality we can assume that every $C^*$- tensor category is a strict rigid $C^*$-tensor category.
\edfn
       
\bxmpl 
The category of finite dimensional Hilbert space, denoted by $Hilb_{f}$, which has finite dimensional Hilbert spaces as the set of objects, morphisms given by linear maps from $\clh$ to $\clk.$ $\clh\otimes \clk$ is the usual tensor product between two Hilbert spaces. Associativity morphism $\alpha_{\clh_{1},\clh_{2},\clh_{3}}$ is the identity map from $(\clh_{1}\ot \clh_{2})\ot \clh_{3}$ to $\clh_{1}\ot (\clh_{2}\ot \clh_{3}).$ $1_{\bbc}$ is unit object of this category $Hilb_{f}.$ $Hilb_{f}$ is a strict rigid $C^{\ast}$ tensor category.
\exmpl
\bxmpl Representation category of a compact group is denoted by $Rep(G).$ Objects of this category are denoted by $(\pi,\clh_{\pi})$, where $\pi:G\to U(\clh_{\pi})$ is a finite dimensional unitary representation of $G$ on $\clh_{\pi}.$ Morphism between two objects $(\pi_{1},\clh_{\pi_{1}})$ $(\pi_{2},\clh_{\pi_{2}})$ is defined by 
   \begin{align*}
    Mor((\pi_{1},\clh_{\pi_{1}}),(\pi_{2},\clh_{\pi_{2}}))=\{T\in\clb(\clh_{\pi_{1}}, \clh_{\pi_{2}}):T\pi_{1}(g)=\pi_{2}T(g),~for~all~g\in G\}.
   \end{align*} Tensor product between two objects $(\pi_{1},\clh_{\pi_{1}})$ $(\pi_{2},\clh_{\pi_{2}})$ is defined by $(\pi_{1}\ot \pi_{2},\clh_{\pi_{1}}\ot \clh_{\pi_{2}}),$  where $\pi_{1}\ot \pi_{2}$ is the usual tensor product of two group representations and $\clh_{\pi_{1}}\otimes \clh_{\pi_{2}}$ is the usual tensor product between two Hilbert spaces. Associativity morphism is the identity map and unit object is the trivial representation on $\bbc.$ $Rep(G)$ is a strict $C^{\ast}$ tensor category.
\exmpl

\bxmpl
Unitary corepresentation category of a CQG denoted by $Corep(\q)$. Objects are given by $(U,\clh_{U}),$ where $U\in B(\clh_{U})\ot\q $ is a finite dimensional unitary corepresentation of $\q$ on $\clh_{U}$. Morphism between two objects $(U,\clh_{U}),(V,\clh_{V})$ is given by \begin{align*}
    Mor((U,\clh_{U}),(V,\clh_{V}))=\{T\in B(\clh_{U},\clh_{V}):(T\ot 1)U=V(T\ot 1)\}
\end{align*} Tensor product of two objects $(U,\clh_{U}),(V,\clh_{V})$ defined by $(U_{13}V_{23},\clh_{U}\ot \clh_{V}).$ Associativity morphism is the identity map and unit object is $Id_{\bbc}\ot 1_{\q}.$ For a projection $P\in End(U,\clh_{U})$, $((P\ot 1)U,P\clh_{U})$ is a subobject of $(U,\clh_{U}).$ $Corep(\q)$ is a strict $C^{\ast}$ tensor category.
\exmpl
\bdfn 
  A fiber functor $F:\clc\to Hilb_{f}$ is a tensor functor which   is   faithful and exact .
\edfn

\bxmpl
Let $F^{Nat}:Corep(\q)\to Hilb_{f}$ is a fiber functor defined by $F^{Nat}(U,\clh_{U})=\clh_{U}$, where $(U,\clh_{U})$ is an object of this category , identity on the morphisms and also $F_{2}(U,V),F_{0}$ are both identity maps.
\exmpl

% We note a standard fact without proof:

%\blmma 
%A linear functor $F$ is faithful if and only if the image of a every simple object is nonzero.
%\elmma

%\begin{proof}
  %  Let assume that, $F$ is faithful. Suppose, for a simple object $U$, $F(U)=0$ that means $End(F(U))=0.$ But we assume that $F$ is faithful and we know $Id_{U}\in End(U).$ It is a contradiction. Therefore, image of a simple object is nonzero.\\
   % suppose, there exist a non zero $f\in Mor(U,V)$ such that $F(f)=0$. $f^{\ast}f\in End(U)$ and $F(f^{\ast}f)=0.$ $End(U)$ is a finite dimensional $C^{\ast}$ that implies it is also a Von Neumann algbera. From this we can conclude that, if we generate a Von Neumann subalgebra by $f^{\ast}f$ then for any projection $P$ of that subalgebra ,$F(P)=0.$ Now, we can choose a  object $V$ of $U$ corresponding to a  projection $P.$ But $F(V)=0$ as $F(p)=0.$ But we assume that image of every simple object is non zero its a contradiction.
%\end{proof}
 %The above lemma is true for any unitary tensor functor because unit object $1$ is a subobject of $F(U\ot\bar{U})\cong F(U)\ot F(\bar{U})$. If for any simple object $U$, $F(U)=0$, then $F(U\ot \bar{U})=0$. So, it is not possible to embed the unit object 1 inside $F(U\ot \bar{U}).$  Any linear functor is exact. Therefore a fiber functor is simply a tensor functor $F:\clc\to Hilb_{f}.$ 
 
% \vspace{2mm}

 \vspace{2mm}
 
 %\vspace{2mm}
 
 \bthm[ Woronowicz's Tannaka-Krein duality\cite{worinv,worole}]
Let $C$ be a a strict rigid $C^{*}$ tensor category and $F$
%:C\rightarrow Hilb_{f}$
be a unitary fiber functor on $C$. Then there exist a CQG $\q$ and unitary monoidal equivalence $E:C\rightarrow Corep(\q)$ such that $F$ is naturally unitary monidally isomorphic to the composition of the canonical fiber functor $Corep(\q)  \rightarrow Hilb_{f}$ with $E$.
\ethm
\subsection{Projective corepresentations and 2-cocycles}
 \bdfn A unitary element $\Omega'\in \q\ot \q,$ is said to be a  right 2-cocycle if it satisfying the equation 
     \begin{align*}
         (id\ot \Delta)(\Omega')(1\ot \Omega')=(\Delta\ot id)(\Omega')(\Omega'\ot 1).
     \end{align*}
     Similarly, a unitary element  $\Omega\in \q\ot \q,$ is said to be a  left 2-cocycle if it satisfying the equation 
     \begin{align*}
         (1\otimes\Omega)(id\otimes \Delta)(\Omega)=(\Omega\otimes 1)(\Delta\otimes id)(\Omega).
     \end{align*}
     \edfn

     \brmrk $\Omega$ is  left 2-cocycle if and only if $\Omega^{\ast}$ is right 2-cocycle.
     \ermrk
\bdfn Let $(\q,\Delta)$ be a $C^{\ast}$-algebraic CQG. A continous left projective coaction of $(\q,\Delta)$ on a Hilbert space $\clh$ consists of a left coaction $\alpha$ of $(\q,\Delta)$ on $K(\clh),$\begin{align*}
    \alpha:K(\clh)\rightarrow M(\q\otimes K(\clh)).
\end{align*}  As it is a non-degenrate $\ast$-homomorphism, we can extend this coaction on $M(K(\clh))=B(\clh).$  Similarly, a continuous right projective coaction of $(\q,\Delta)$ on a Hilbert space $\clh$ consists of a right coaction $\beta$ of $(\q,\Delta)$ on $B(\clh)$.\edfn
\bdfn Let $\delta$ be a continuous left (right) projective coaction. We say that $\delta$ is \emph{cleft} if there exists a unitary $U\in \q \ot B(\clh)$  ($B(\clh)\ot\q$) such that  $\delta(a) = u^{\ast}(1\otimes a)u$  ($\delta(a)=u(a\otimes 1)u^{\ast}$).
 \vspace{1mm}
 A unitary $U\in M(K(\clh)\otimes \q)$ is said to be a continuous right projective  corepresentation if $(Id\otimes \Delta_{\Omega^{\ast}})(U)=U_{12}U_{13}$ for a right 2-cocycle $\Omega^{\ast}\in \q\ot \q$ and a continuous left projective $\Omega$ corepresentation if  $(Id\otimes _\Omega\Delta)(U)=U_{12}U_{13}$ for a left 2-cocycle $\Omega\in (\q\ot\q).$
 \edfn
 \bppsn (proposition 3.1.9 ~ of~\cite{proj2DeMaNe}) Let $\delta: B(\clh) \rightarrow B(\clh) \ot \q $ be a cleft right projective coaction. Then  there exists a right $2$-cocycle $\Omega \in \q\ot \q$ and a unitary  right projective  corepresentation  $U\in B(\clh) \ot \q$ such that $\delta(a) = U(a\otimes 1)U^*.$
 \eppsn
 %\brmrk If $\delta:B(\clh)\rightarrow L^{\infty}(\q)\bar\otimes B(\clh) $ be a cleft left measurable projective corepresentation if and only if there is a right 2-cocycle $\Omega^{\ast}\in L^{\infty}(\q)\bar\otimes L^{\infty}(\q)$ and a measurable left projective corepresentation $U\in B(H)\bar\otimes L^{\infty}(\q)$ such that $\delta(a)=\sigma( U^{\ast}(a\otimes 1)U),$ where $\sigma$ is the flip map from $B(\clh)\bar{\ot} L^{\infty}(\q) ~ to~ L^{\infty}(\q)\bar{\ot}B(\clh).$
 %\ermrk
 \bppsn (Theorem 3.1.12 of \cite{proj2DeMaNe}) Any continuous  left/right projective coaction of $(\q,\Delta)$ is cleft, which means any  right projective coaction $\alpha$  corresponds to a unitary right projective  corepresentation $U\in M(\clk(H)\otimes\q)$ such that $\alpha=Ad_{U}$ and  any left projective coaction $\beta$ corresponds to a left projective corepresentation $V\in M(\clk(H)\otimes\q)$ such that $\beta=\sigma\circ Ad_{V^{\ast}},$ where $\sigma$ is the flip map.
 \eppsn
% \brmrk Note the slight confusion in the terminology as a projective left/right unitary corepresentation can mean either a homomorphism $\delta$  or a corresponding unitary. We often call a unitary to be a corepresentation without explicitly mentioning the cocycle $\Omega$ if it is understood from the context. However, we will mostly use unitary picture of projective corepresentation. For a classical compact group $G,$  left/right projective representation are the same thing. However, for a general CQG this is not so. Thus, it is natural to consider special class left/right/biprojective projective corepresentations.
 %\ermrk
 We now begin to initiate some of the special classes of projective corepresentations.
 \bdfn 
   We say a unitary $U\in M(K(\clh)\otimes \q)$ is  biprojective  corepresentation if it is a both left and right projective corepresentation (possibly with different 2-cocycles) of $\q.$
 \edfn
 
 \bdfn Let $U\in M(K(\clh_{u})\otimes \q),V\in M(K(\clh_{v})\otimes \q)$ be two right/left projrctive corepresentations of $\q.$ Morphism between $U,V$ are given by \begin{align*}
     Mor(U,V)=\{T\in B(\clh_{u},\clh_{v}):(T\otimes 1)U=V(T\otimes 1)\}.
 \end{align*} $(U,\clh)$ is called irreducible if $Mor(U,V)=\bbc Id_{\clh_{U}}.$
 \edfn
\blmma(Lemma (3.2.5) of \cite{proj2DeMaNe}) If $(U, \clh_{u})$ and $(V, \clh_{v})$ are two irreducible  left/right projective corepresentations of $\q$, then either $U$ is not unitary equivalent to $V$ and $Mor(U,V)=(0)$; or $U$ is unitary equivalent to $V$ and $Mor(U,V)$ is a $1$-dimensional subspace of $B(\clh_{u}, \clh_{v})$.\elmma
\blmma\label{dirs}(Lemma (3.2.6) 0f\cite{proj2DeMaNe})  Every  left/right projective corepresentation $U$ of $\q$ decomposes into a direct sum of irreducible left/right projective corepresentations.\elmma
\bppsn Let $\q$ be a compact quantum group and $\Omega$ a left $2$-cocycle on $\q$. Defining $V^{\Omega} =\Omega V_{\q}$, where $V_{\q}$ is right regular corepresentation of $\q$, the following properties hold:
			\begin{enumerate}
				\item[i)] For all $x\in \q $ and $\xi\in H_{u}$ we have $V^{\Omega}(\Lambda(x)\otimes \xi)={}_{\Omega}\Delta(x)(\xi_{\q}\otimes \xi)$.
				\item[ii)]For all $x\in \q $ we have ${}_{\Omega}\Delta(x)=V^\Omega(x\otimes 1)V_{\q}^*$.
				\item[iii)] The following identity holds: $(id\otimes {}_{\Omega}\Delta)(V^\Omega)=V^\Omega_{12}V^\Omega_{13}$, so $V^{\Omega} \in \mathcal{B}(L^2(\q))\otimes \q$ is an $\Omega$-representation.
				\item[iv)] The following pentagonal equation holds: $V^\Omega_{12}V^\Omega_{13}(V_{\q})_{23}=V^{\Omega}_{23}V^\Omega_{12}$.
			\end{enumerate}
			The unitary $V^\Omega$ is a \emph{left regular projective corepresentation of $\q$ on $L^2(\q)$ with respect to $\Omega$} or  \emph{left regular $\Omega$ corepresentation of $\q$} on $L^2(\q)$. \eppsn
          	\brmrk Similarly, defining $W^{\Omega} = W_{\q}\Omega^*$, we have that  $(W^{\Omega})^*(\xi\otimes \Lambda(x))={}_{\Omega}\Delta(x)(\xi\otimes \xi_{\q})\mbox{,}$
		for all $x\in \q$ and $\xi\in L^2(\q)$. For all $x\in \q$, we have ${}_{\Omega}\Delta(x)=(W^\Omega)^*(1\otimes x)W_{\mathbb{G}}$ and the pentagonal equation: $(W_{\mathbb{G}})_{12}W^\Omega_{13}W^\Omega_{23}=W^{\Omega}_{23}W^\Omega_{12}$ and the following identity holds: $(\Delta_{\Omega^*}\otimes id)(W^\Omega)=W^\Omega_{13}W^\Omega_{23}$, so $\Sigma W^{\Omega}\Sigma$ is an $\Omega^*$-projective corepresentation.
					
		The unitary $W^\Omega$ is called \emph{right  regular projective $\Omega$ corepresentation of $\q$ on $L^2(\q)$ with respect to $\Omega$} or simply \emph{right regular $\Omega$ corepresentation of $\q$ on $L^2(\q)$}. \ermrk
        \bppsn[Twisted Peter-Weyl theorem I]\label{theo.PeterWeyl1Meas}(Theorem (3.2.12) of \cite{proj2DeMaNe})
		Let $\q$ be a CQG and $\Omega$ a $2$-cocycle. The right regular projective  corepresentation $(V^\Omega, L^2(\q)$ contains all irreducible $\Omega$-representations of $\q$ in its direct sum decomposition.\eppsn
\bppsn 
     [Twisted Schur's orthogonality relations]\label{theo.TwistedOrthogonalityRel}(Theorem (3.2.13 of \cite{proj2DeMaNe})
		Let $\q$ be a CQG and $\Omega$ a left $2$-cocycle on $\q$. Let $\{u^x\}_{x\in Irr(\q, \Omega)}$ be a complete set of mutually inequivalent, irreducible left projective $\Omega$-corepresentations, with fixed bases for the associated Hilbert spaces $\clh_x$. For each $x\in Irr(\mathbb{G}, \Omega)$ there exists a positive trace class operator $F^x\in\mathcal{B}(\clh_x)$ with zero kernel such that the following orthogonality relations hold:
		$$h_{\q}\big((u^y_{kl})^*u^x_{ij}\big)=\delta_{x y}\delta_{l j}F^x_{ik},$$
		for every $x,y\in Irr(\q, \Omega)$, $i,j=1,\ldots, n_x$ and $k,l=1,\ldots, n_y$, $h_{\q}$ is the Haar state of $(\q,\Delta)$\eppsn 
        \bthm [Twisted Peter-Weyl theorem II]\label{theo.TwistedPeterWeyL2}(Theorem (3.2.14) of \cite{proj2DeMaNe})
		Let $\q$ be a CQG and $\Omega$ a left $2$-cocycle on $\q$. We have a unitary transformation $L^2(\q)\cong \underset{x\in Irr(\q, \Omega)}{\bigoplus} \clh_x\otimes \overline{\clh_{x}}$ such that $\Lambda(u^x_{ij})\mapsto \sqrt{F^x_i}\ \xi^x_{i}\otimes \overline{\xi^{x}_{j}}$, for all $j=1,\ldots, n_x$, $x\in \text{Irr}(\q, \Omega)$.
        Similar statement hold for  right projective corepresentation of $\q$ for a corresponding right 2-cocycle $\Omega^{\ast}.$
        \ethm
        \bppsn 
	[Twisted Maschke's theorem]\label{theo.AverageInterwinersCont}
		Let $\q$ be a compact quantum group and $\Omega$ a  $2$-cocycle on $\q$. Let $(U, \clh_{u})$ and $(V, \clh_{v})$ be two right projective  $\Omega$-corepresentations of $\q$.
		\begin{enumerate}
\item[i)] If $T: \clh_{u}\longrightarrow \clh_{v} $ is a linear compact operator, then the average intertwiner $T^{'}$ with respect to $T$ is again compact. 
			\item[ii)] The $C^*$-algebra $\mathcal{D}_u = \mathcal{K}(\clh_{u})^{\delta_u}$ acts non-degenerately on $\clh_{u}$, that is, $[\mathcal{D}_{u}\clh_{u}]=\clh_{u}$.
\item[iii)]\label{It.DecompIrr} If $U$ is an irreducible object, then $U$ is finite dimensional. 
 
		\end{enumerate}\eppsn

\section{Unitarity of the conjugate $U^{c}$ for a projective corepresentation $U$}

       The key step of this construction is the fact that for any unitary projective corepresentation $U$ on a finite-dimensional Hilbert space, $U^{c}$      is unitary for a suitable choice of inner product. This will be proved in this section.
        It has already been seen in lemma \ref{dirs}  that any continuous projective corepresentation can be expressed as a direct sum of irreducible projective corepresentations and  any irreducible projective corepresentation is always finite dimensional. Therefore it is natural to study  finite-dimensional projective corepresentations. From now on. we will assume that any projective corepresentation considered by us is  finite dimensional.  
        %So ,from now on we work on a finite-dimensional Hilbert space.\\

        Fix an orthonormal basis $\{e_{i}:i=\mathrm {1,..., n}\}$ of a finite-dimensional Hilbert space $\clh$. Hence $\{\bar{e_{i}}:i=1,..., n\}$ is an orthonormal basis for the Hilbert space $\bar{\clh}$, where the inner product $< , >_{\overline {\clh}}$  defined by
     $<\bar{a},\bar{b}  >_{\overline {\clh}}$ =$<b,a>_{\clh}$ .  
     
     Let $\phi:B(\clh) \rightarrow \clh \otimes \overline{\clh}$ be the linear map defined by $\phi(e_{ij})=e_{i}\otimes\overline{e_{j}}$, where $e_{ij}$ are the matrix units of $B(\clh)$, given by $e_{ij}.$ Similarly, let $\overline{e_{ij}}$ be the matrix units of $B(\bar{\clh})$ given by $\bar{e_{ij}}(\bar{e_{k}})=\bar{e_{i}}\delta_{kj}.$  Clearly $\phi$ is a vector space isomorphism.

     Let $U\in B(\clh)\ot\q$ be a unitary element. Define $U^{c}=(j\ot id)(U^{\ast}),$ where $j:B(\clh)\rightarrow B(\bar{\clh})$ is the canonical complex conjugation map. $(U\ot U^{c})=U_{13}U^{c}_{23}\in B(\clh\ot \bar{\clh})\ot\q\subseteq B(\clh\ot\bar{\clh}\ot L^{2}(\q)),$ where $L^{2}(\q)$ is the GNS space for the Haar state (say $h$) on $\q.$
\blmma\label{tr} 
    \begin{align*}
        Ad_{U}(a)=(\phi^{-1}\otimes id_{\q})(U\otimes U^{c})(\phi\otimes Id_{\q})(a\otimes 1)
    \end{align*}
where $a\in B(H)$ and 1 is the identity element of $\q$ viewed as a cyclic vector of $L^{2}(\q)$.
\elmma
\begin{proof}
    Assume that $U=\sum e_{i,j}\otimes u_{i,j}$,where $e_{ij}$ are the matrix units of $B(\clh)$ and $u_{i,j}\in \q$.
For $a=e_{i_{0},j_{0}}$,
\begin{align*}
Ad_{U}(e_{i_{0}j_{0}})=&U(e_{i_{0}j_{0}}\otimes 1)U^{*}\\
=&\sum_{i_1,j_{1},i_{2},j_{2}}(e_{i_{1}j_{1}} \otimes u_{i_{1}j_{1}})(e_{i_{0}j_{0}}\otimes 1)(e_{j_{2}i_{2}}\otimes u^{*}_{i_{2}j_{2}})\\
=&\sum e_{i_{1}j_{1}}e_{i_{0}j_{0}}e_{j_{2}i_{2}}\otimes u_{i_{1}j_{1}}u^{*}_{i_{2}j_{2}}\\
=&\sum_{i_{1},i_{2}} e_{i_{1}i_{2}}\otimes u_{i_{1}i_{0}}u^{*}_{i_{2}j_{0}}.
\end{align*}
We know $U\otimes U^{c}=U_{13}U^{c}_{23}=\sum e_{i_{1}j_{1}}\otimes \bar{e}_{i_{2}j_{2}}\otimes u_{i_{1}j_{1}}u^{*}_{i_{2}j_{2}}$.

Now,\begin{align*}
&(\phi^{-1}\otimes Id_{\q})(U\otimes U^{c})(\phi \otimes Id_{\q})(e_{i_{0}j_{0}}\otimes 1)\\
=&(\phi^{-1}\otimes Id_{\q})(U\otimes U^{c})(\phi(e_{i_{0}j_{0}}) \otimes 1)\\
=&(\phi^{-1}\otimes Id_{\q})(U\otimes U^{c})(e_{i_{0}}\otimes \bar{e}_{j_{0}}\otimes 1)\\
=&(\phi^{-1}\otimes Id_{\q})(\sum e_{i_{1}j_{1}}(e_{i_{0}})\otimes \bar{e}_{i_{2}j_{2}}(\bar{e}_{j_{0}})\otimes u_{i_{1}j_{1}}u^{*}_{i_{2}j_{2}})\\
=&(\phi^{-1}\otimes Id_{\q})(\sum e_{i_{1}}\delta_{i_{0}j_{1}}\otimes \bar{e_{i_2}}\delta_{j_{0}j_{2}}\otimes u_{i_{1}j_{1}}u^{*}_{i_{2}j_{2}})\\
=&(\phi^{-1}\otimes Id_{\q})(\sum e_{i_{1}}\otimes \bar{e_{i_{2}}}\otimes u_{i_{1}i_{0}}u^{*}_{i_{2}j_{0}})\\
=&\sum \phi^{-1}(e_{i_{1}}\otimes \bar {e_{i_{2}}})\otimes u_{i_{1}i_{0}}u^{*}_{i_{2}j_{0}}\\
=&\sum e_{i_{1}i_{2}}\otimes u_{i_{1}i_{0}}u^{*}_{i_{2}j_{0}}.
\end{align*}
As the matrix units $e_{i_{0}j_{0}}$ form a basis of $B(\clh)$ , it follows that $Ad_{U}(a)=(\phi^{-1}\otimes Id_{\q})(u\otimes u^{c})(\phi\otimes Id_{\q})(a\otimes 1)$ for any $a\in B(\clh)$.
\end{proof}

Let $\delta$ be a right coaction of $\q$ on $B(\clh)$ that corresponds to a right projective corepresentation $U\in B(\clh)\ot\q.$ We can write $\delta=Ad_{U}.$

Let $\psi:B(\clh)\rightarrow \bbc $ be the linear function defined by $\psi(a)=(Tr\otimes h)\delta(a),$ where $h$ is the Haar state on $\q$ as mentioned before.
%$\psi$ is a faithful state on $B(H)$, hence there exists a positive invertible matrix $R$ in $B(H)$ s.t. $\psi(a)=Tr(Ra)$, where $a$ is in $B(H)$.
\blmma\label{Tr}
 We have
\begin{align*}
    (\psi\otimes Id_{\q})\delta(a)=\psi(a)1_{\q},
\end{align*} where $a\in B(\clh).$ 
\elmma

\begin{proof}
For $a\in B(\clh)$, we have

    \begin{align*}
&(\psi\otimes Id_{\q})\delta(a)\\
=&((Tr\otimes h)\delta \otimes Id_{\q})\delta (a)\\
=&(Tr\otimes h\otimes Id_{\q})(\delta \otimes Id_{\q})\delta (a)\\
=&(Tr\otimes h\otimes Id_{\q})(Id_{B(H)}\otimes \Delta)\delta(a)\\
=&(Tr\otimes (h\otimes Id_{\q})\Delta)\delta(a)\\
=&(Tr\otimes h.1_{\q})\delta(a)\\
=&(Tr\otimes h)\delta(a)1_{\q}\\
=&\psi(a)1_{\q}.
\end{align*}
%So, $(\psi\otimes Id_{\q})\delta(a)=\psi(a)1_{\q}$.
\end{proof}
\brmrk
As $\psi$ is a faithful state on $B(\clh)$ there exists a positive invertible linear operator $R\in B(\clh)$ such that $\psi(a)=Tr(Ra).$
Now, we define an inner product $<.,.>_{R}$ on $\bar{\clh}$ st. $<\bar{a},\bar{b}>_{R}=<\overline{R^{1/2}(a)},\overline{R^{1/2}(b)}>$\label{ip}. Let us denote this Hilbert space by $\bar{\clh}_{R}$ and %we can give a new inner product on $B(\clh)$ by using this positive operator $R$ and for any two operators
define a  new inner product on $B(\clh)$ by:
$<S,T>_{new}=\psi(T^{*}S)=Tr(RT^{*}S)$. It is easy to check that under this new inner product ($B(\clh),<S,T>_{new}$) is isometrically isomorphic to $\clh\otimes \bar{\clh}_{R}.$
\ermrk
\bthm\label{cop}
Define the linear map
$\delta_{U}:\clh\otimes \bar{\clh}_{R}\rightarrow \clh\otimes\bar{\clh}_{R}\otimes \q $ given by
\begin{align*}
    \delta_{U}=(U\otimes U^{c}).
\end{align*}
Then $\delta_{U}$ is a unitary corepresentation of the CQG $(\q,\Delta).$
\ethm
\begin{proof}
  As $R$ is a positive matrix in $B(\clh)$ hence $R$ is  diagonalizable. Without loss of generality, we assume that $e_{i}$ are orthonormal eigenvectors of $R$. Therefore, $\{e_{i}\otimes \bar{e_{j}}:i,j ={1\cdots n}\} $ is an orthonormal basis for $ H\otimes \bar{H_{R}}$.

  We will prove that $<\delta_{U} (e_{i_{0}}\otimes \bar{e_{j_{0}}}),\delta_{U} (e_{l_{0}}\otimes \bar{e_{m_{0}}})>=<e_{i_{0}}\otimes \bar{e_{j_{0}}},e_{l_{0}}\otimes \bar{e_{m_{0}}}>.$
\begin{align*}
  &<\delta_{U}(e_{i_{0}}\otimes \bar{e_{j_{0}}}),\delta_{U}(e_{l_{0}} \otimes \bar{e_{m_{0}}} )>\\
  =&<(U\otimes U^{c})(e_{i_{0}}\otimes \bar{e_{j_{0}}}\otimes 1),(U\otimes U^{c})(e_{l_{0}}\otimes \bar{e_{m_{0}}}\otimes 1)\\
  =&<\sum_{i_{1},i_{2}}e_{i_{1}}\otimes \bar{e_{i_{2}}}\otimes u_{i_{1}i_{0}}u^{*}_{i_{2}j_{0}}, \sum_{l_{1},l_{2}}e_{l_{1}}\otimes \bar{e_{l_{2}}}\otimes u_{l_{1}l_{0}}u^{*}_{l_{2}m_{0}}>\\
  =&\sum_{i_{1},i_{2},l_{1},l_{2}}<e_{i_{1}}\otimes \bar{e_{i_{2}}},e_{l_{1}}\otimes \bar{e_{l_{2}}}>(u_{i_{1}i_{0}}u^{*}_{i_{2}j_{0}})^{*}u_{l_{1}l_{0}}u^{*}_{l_{2}m_{0}}\\
  =&\sum<e_{i_{1}},e_{l_{1}}><\bar{e_{i_{2}}},\bar{e_{l_{2}}}>_{\bar{H_{R}}}u_{i_{2}j_{0}}u^{*}_{i_{1}i_{0}}u_{l_{1}l_{0}}u^{*}_{l_{2}m_{0}}\\
  =&\sum \delta_{i_{1}l_{1}}<Re_{l_{2}},e_{i_{2}}>u_{i_{2}j_{0}}u^{*}_{i_{1}i_{0}}u_{l_{1}l_{0}}u^{*}_{l_{2}m_{0}}\\
  =&\sum_{j_{2},l_{2}}<Re_{l_{2}},e_{i_{2}}>u_{i_{2}j_{0}}(\sum u^{*}_{i_{1}i_{0}}u_{i_{1}l_{0}})u^{*}_{l_{2}m_{0}}\\
  =&\sum \delta_{l_{0}i_{0}}<Re_{l_{2}},e_{i_{2}}>u_{i_{2}j_{0}}u^{*}_{l_{2}m_{0}}.
\end{align*}
From Lemma (\ref{Tr}), we know that $(\psi\otimes Id)\delta(a)=\psi(a).$\\
For $a=e_{m_{0}j_{0}},$ we have

\begin{align*}
    &(\psi\otimes Id)\delta(e_{m_{0}j_{0}})\\
    =&(\psi\otimes Id)Ad_{U}(e_{m_{0}j_{0}})\\
    =&(\psi\otimes Id)(\sum e_{i_{1}i_{2}}\otimes u_{i_{1}m_{0}}u^{*}_{i_{2}j_{0}})\\
    =&\sum \psi(e_{i_{1}i_{2}})u_{i_{1}m_{0}}u^{*}_{i_{2}j_{0}}\\
 =&\sum<Re_{i_{1}},e_{i_{2}}>u_{i_{1}m_{0}}u^{*}_{i_{2}j_{0}}=\psi(e_{m_{0}j_{0}})1_{\q}.
\end{align*}
Taking adjoint on both side of the equation we will get, $\sum<Re_{i_{1}},e_{i_{2}}>u_{i_{2}j_{0}}u^{*}_{i_{1}m_{0}}=\psi(e_{m_{0}j_{0}})1_{\q}$.\\
It follows that,
\begin{align*}
    &\sum\delta_{l_{0}i_{0}}<Re_{l_{2}},e_{i_{2}}>u_{i_{2}j_{0}}u^{*}_{l_{2}m_{0}}\\
&=\delta_{l_{0}i_{0}}\psi(e_{m_{0}j_{0}})1_{\q}\\
&=<e_{i_{0}},e_{l_{0}}><R^{1/2}e_{m_{0}},R^{1/2}e_{j_{0}}>1_{\q}\\
&=<e_{i_{0}},e_{l_{0}}><\bar{e_{j_{0}}},\bar{e_{m_{0}}}>_{\bar{H_{R}}}1_{\q}\\
&=<e_{i_{0}}\otimes \bar{e_{j_{0}}},e_{l_{0}}\otimes \bar{e_{m_{0}}}>_{H\otimes \bar{H_{R}}}1_{\q}.
\end{align*}
Thus $<\delta (e_{i_{0}}\otimes \bar{e_{j_{0}}}),\delta (e_{l_{0}}\otimes \bar{e_{m_{0}}})>=<e_{i_{0}}\otimes \bar{e_{j_{0}}},e_{l_{0}}\otimes \bar{e_{m_{0}}}>$ and using the fact that matrix units form a basis of $\clh\ot\bar{\clh}_{R}$ we get $<\delta(a),\delta(b)>=<a,b>1_{\q}$, for any $a,b\in \clh\otimes \bar{\clh}_{R}$.

 Now, we will prove $(\delta_{U}\otimes id)\delta_{U}=(Id\otimes \Delta)\delta_{U}$.\\
 we know from Lemma (\ref{tr}),
  $Ad_{U}(\phi^{-1}(a))=(\phi^{-1}\otimes Id)(U\otimes U^{c})(a\otimes 1)$, where $a\in \clh\otimes \bar{\clh},$ which gives:
 \begin{align*}
(Ad_{U}\otimes Id)Ad_{U}(\phi^{-1}(a))=&(id\otimes \Delta)Ad_{U}(\phi^{-1}(a))\\
=&(id\otimes \Delta)(\phi^{-1}\otimes Id_{\q})\delta_{U}(a)\\
=&(\phi^{-1}\otimes Id_{\q}\otimes Id_{\q})(Id\otimes \Delta)\delta_{U}(a).\\     
\end{align*}
Then, it easily follows that 
\begin{align*}
   &(Id\otimes \Delta)\delta_{U}(a)\\
   =& (\phi\otimes Id_{\q}\otimes Id_{\q})(Ad_{U}\otimes Id)Ad_{U}(\phi^{-1}(a))\\
   =&(\phi\otimes Id_{\q}\otimes Id_{\q})(Ad_{U}\otimes Id)(\phi^{-1}\otimes id)\delta_{U}(a)\\
    =&((\phi\otimes Id_{\q})\otimes Id_{\q})(Ad_{U}\phi^{-1}\otimes Id)\delta_{U}(a)\\
   =&((\phi\otimes Id_{\q})\otimes Id_{\q})((\phi^{-1}\otimes Id)\delta_{U}\otimes Id_{\q})\delta_{U}(a)
    =(\delta_{U}\otimes Id_{\q})\delta_{U}(a).\\
   \end{align*}

   Note that $\delta=Ad_{U}$ is a coaction of CQG $(\q,\Delta)$ on the finite-dimensional space $B(\clh),$ hence from the general theory of CQG coaction, span$\{Ad_{U}(x)(1\ot q):x\in B(\clh),q\in \q\}=B(\clh)\ot \q.$  As $\phi^{-1}$ is a vector space isomorphism, we have span$\{Ad_{U}(\phi^{-1})(a)(1\ot q):a\in \clh\ot\bar{\clh}_{R},q\in \q\}=B(\clh)\ot \q.$ This proves that $\delta_{U}(\clh\ot\bar{\clh})(1\ot\q)$ is total in $\clh\ot\bar{\clh}_{R}\ot\q.$ Hence, $\delta_{U}$ is a unitary corepresentation.
% We know that span\{$Ad_{U}\phi^{-1}(a).\q:\forall a\in H\otimes \bar{H_{R}}$\} is dense in $B(H)\otimes \q$ and $\phi$ is a vector space isomorphism from $B(H)$ to $H\otimes\bar{H}.$ So,from those two facts,We can conclude that $\delta(H\otimes \bar{H_{R}}).\q$ is dense in $H\otimes \bar{H_{R}}\otimes \q.$
 %$\delta$ is a unitary corepresentation.  
\end{proof}
   
 \bcrlre
     $U^{c}$ is a unitary element of $B(\bar{\clh}_{R})\otimes \q$.
 \ecrlre
 
     % $\delta_{u}$ is a unitary corepresentation of $(\q,\Delta)$ on the Hilbert space $H\otimes \bar{H_{R}}$. So, there exists a unitary \begin{align*}
         %& X:H\otimes \bar{H_{R}}\otimes \q\rightarrow H\otimes \bar{H_{R}}\otimes \q \\
       %  &\textrm{such that } X(a\otimes \bar{b}\otimes c)=\delta(a\otimes\bar{b}).c
      %\end{align*}
% , for all $a\in H,b\in \bar{H}_{R}, c\in \q$ and also $X\in M(K(H\otimes \bar{H_{R}})\otimes \q)=B(H)\otimes B(\bar{H_{R}})\otimes \q$.
 %$X$ is a unitary element of the $C^{*}$ algebra $B(H)\otimes B(\bar{H_{R}})\otimes \q$ st.$(Id\otimes \Delta )(X)=X_{12}X_{13}$. One can easily check that $X(a\otimes b\otimes c)=(u\otimes u^{c})(a\otimes b\otimes c)$. So,$u\otimes u^{c}$ is a unitary element of the $C^{*}$ algebra $B(H)\otimes B(\bar{H_{R}})\otimes \q$ and we know that $u_{13}$ is a unitary element. Hence, $u_{13}^{*}(u\otimes u^{c})=u_{23}^{c}$ is a unitary element of the $C^{*}$ algebra  $B(H)\otimes B(\bar{H_{R}})\otimes \q$ . Now,we can conclude that $u^{c}$ is a unitary element of $B(\bar{H_{R}})\otimes \q.$
\begin{proof}
This follows from the fact that 
%$U_{13}^{*}(U\otimes U^{c})=U_{23}^{c}$ and $U^{c}_{23}$ is a unitary. Hence,
$(1\ot U^{c})=(U^{-1}\ot 1)\delta_{U}$ which is unitary, as both $U^{-1}$ and $\delta_{U}$ are unitary on $\clh\ot \bar{\clh}_{R},$ which implies $U^c$ is unitary on $\bar{\clh}_{R}$ as claimed.
 
%U_{12}^{-1}U_{13}U^{c}_{23},$ which is unitary.
 \end{proof}

Now, we define a linear map $T:\bar{\clh}\rightarrow \bar{\clh}$ such that $T(\bar{x})=\overline{R^{(-1/2)}(x)}$. If both the domain and range have the same old inner product structure, then $T$ is a positive invertible map.
If the domain space has the original inner product structure of $\bar{\clh}$ and the range space has the inner product $< , >_{{\bar{\clh}}_{R}}$ as in remark (\ref{ip}), which is induced from the positive matrix $R$, then $T$ is an isometry from $\clh$ to $\bar{\clh}_{R}$ because $<T(\bar{x}),\overline{y}>_{{\bar{\clh}}_{R}}=<\overline{R^{-(1/2)}(x)}, 
 \overline{y}>_{{\bar{\clh}}_{R}}=<\overline{x}, \overline{R^{1/2}(y)}>$  implies that $T^{*}(y)=\overline{R^{1/2}(y)}$ and $T^{*}T=Id_{\bar{H}}$. $T$ is an isometry from $(\bar{\clh},< , >) $ to $(\bar{\clh},< , >_{{\bar{\clh}}_{R}})$ and $TT^{*}=Id_{\bar{\clh}}$.

\blmma
    $(T^{-1}\otimes Id_{\q})U^{c}(T\otimes Id_{\q})$ is a unitary element of $B(\bar\clh)\otimes \q$, where the inner product structure in $\bar{\clh}$ is the natural inner product structure.
\elmma
\begin{proof}
    $U^{c}$ is a unitary element in $B(\bar{\clh}_{R})\otimes \q$. $T$ is an isometry from $\bar{\clh}$ to $\bar{\clh}_{R}$ and an invertible map. By using these three facts, we can conclude that $(T^{-1}\otimes Id_{\q})U^{c}(T\otimes Id_{\q})$ is a unitary element.
\end{proof}

 We assume that $T^{-1}=\rho^{1/2}$, where $\rho(\bar{x})=\overline{R(x)}.$
\blmma\label{rightp}
    If $U$ is a unitary right projective  corepresentation of $\q$ on $\clh$, then $(\rho^{1/2}\otimes 1)U^{c}(\rho^{-1/2}\otimes 1)$ is a unitary left projective corepresentation on $\bar{\clh}.$
\elmma
\begin{proof}
    If $U$ is a unitary right projective  corepresentation with the associated right 2-cocycle $\Omega$ then it follows that $(Id\otimes \Delta_{\Omega^{*}})(U)=U_{12}U_{13}$. Now, one can easily conclude that $(Id\otimes {_\Omega}\Delta)(U^{*})=U^{*}_{13}U^{*}_{12}.$\\
We know $U^{c}=(j\otimes Id)(U^{*})$, where $j$ is the conjugation map from $B(\clh)$ to $B(\bar{\clh)})$.\\
This implies, \begin{align*}
   (Id\otimes {_\Omega}\Delta)(U^{c})=& (Id\otimes {_\Omega}\Delta)(j\otimes Id)(U^{*})\\
   =&(j\otimes Id\otimes Id)(Id\otimes {_\Omega}\Delta)(U^{*})\\
   =&(j\otimes Id\otimes Id)(U^{*}_{13}U^{*}_{12})\\
   =&(j\otimes Id\otimes Id)(U^{*}_{12})(j\otimes Id\otimes Id)(U^{*}_{13})\\
   =&U^{c}_{12}U^{c}_{13}.
   \end{align*}
   From this we conclude that  
   \begin{align*}
     (Id\otimes {_\Omega}\Delta)((p^{1/2}\otimes 1)u^{c}(p^{-1/2}\otimes 1))=&(p^{1/2}\otimes 1\otimes 1)u^{c}_{12}u^{c}_{13}(p^{-1/2}\otimes 1\otimes 1)\\
     =&(p^{1/2}\otimes 1\otimes 1)u^{c}_{12}(p^{-1/2}\otimes 1\otimes 1)(p^{1/2}\otimes 1\otimes 1)u^{c}_{13}(p^{-1/2}\otimes 1\otimes 1)\\
     =&((\rho^{1/2}\otimes 1)u^{c}(\rho^{-1/2}\otimes 1))_{12}((\rho^{1/2}\otimes 1)u^{c}(\rho^{-1/2}\otimes 1))_{13},
 \end{align*}
where $(\rho^{1/2}\otimes 1)u^{c}(\rho^{-1/2}\otimes 1)$ is a unitary left projective  corepresentation on $\bar{\clh}$.
\end{proof}

%\bdfn A finite dimensional unitary $\Omega^{*}$ representation $u\in B(H)\otimes \q$ is said to be the right projective corepresentation  of $(\q,\Delta)$ and a finite-dimensional unitary representation $\Omega^{'}$ $u\in B(H)\otimes \q$ is said to be the left projective representation  of $(\q,\Delta)$. It is said to be biprojective if it is a right projective representation and $\bar{u_{r}}=(\rho^{1/2}\otimes 1)u^{c}(\rho^{-1/2}\otimes 1)$ is a right projective representation. Here, $\Omega,\Omega^{'}$ are both 2-cocycles of $\q$.\edfn

\blmma \label{cor}
     Let $U\in B(H)\ot \q$ is a unitary.Then, followings are equivalent: 
     \begin{enumerate}
         \item [1)]If $U$ is a right projective corepresentation.  \item[2)]There exists a positive matrix $\rho\in B(\bar{H})$ st.  $U\otimes ((\rho^{1/2}\otimes 1)U^{c}(\rho^{-1/2}\otimes 1)) $  is a unitary corepresentation of $(\q,\Delta).$
         \item[3)]$U\ot U^{c}$ is an invertible corepresentation on $H\ot \bar{H}.$
     \end{enumerate} 
\elmma
\begin{proof}
    $(1\Rightarrow 2) $ 
    
     Assume that $U$ is a right projective  corepresentation with corresponding 2-cocycle $\Omega.$ From the lemma (\ref{rightp}) , we know that there exists a positive matrix $\rho$ such that $\bar{U_{r}}:=(\rho^{1/2}\otimes 1)U^{c}(\rho^{-1/2}\otimes 1)$ is a left projective corepresentation.
Now,
\begin{align*}
    (Id\otimes \Delta)(U\ot \bar{U_{r}})=&(Id\otimes \Delta)(U_{13}\bar{U_{r}}_{23})\\
    =&(Id\otimes \Delta)(U_{13})(Id\otimes \Delta)(\bar{U_{r}}_{23})\\
    =&U_{13}U_{14}\Omega .\Omega^{*}\bar{U_{r}}_{23}\bar{U_{r}}_{24}\\
    =&U_{13}\bar{U_{r}}_{23}U_{14}\bar{U_{r}}_{24}\\
    =&(U\ot \bar{U_{r}})_{12}(U\ot \bar{U_{r}})_{13}.
\end{align*}
 This proves that $(U\ot \bar{U_{r}})$ is a unitary  corepresentation of $(\q,\Delta)$.

\vspace{2mm}

$(2\Rightarrow 3)$

$(U\ot U^{c})=(1\ot (\rho^{1/2}\ot 1_{\q})(U\ot \bar{U_{r}})(1\ot (\rho^{-1/2}\ot 1_{\q}).$ As $U\ot \bar{U}_{r}$  is a unitary corepresentation therefore $U\ot U^{c}$ is an invertible corepresentation of $(\q,\Delta)$.

\vspace{2mm}

$(3\Rightarrow 1)$

 Let $U\ot U^{c}$ be an invertible corepresentation of $\q.$ As we observe, $Ad_{U}=(\phi^{-1}\otimes 1_{\q})(U\otimes U^{c})(\phi\otimes 1_{\q})$( Proposition \ref{tr}). Now, the proof of coassociativity in Theorem \ref{cop} applies verbatim in this case as well. Hence $Ad_{U}$ is a coaction of $\q$ on $B(\clh).$ $U$ is a right projective  corepresentation corresponding to a 2-cocycle $\Omega$ of $\q.$\end{proof}
  % Let $U\ot U^{c}$ is an invertible corepresentation of $\q$.
%Now, we define a map  $Ad_{U}:B(H)\rightarrow B(H)\otimes \q$ such that $Ad_{U}(a)=u(a\otimes 1)u^{*}$, for all $a\in B(H)$.
%If $u=\sum e_{ij}\otimes u_{ij}$ and $(u\otimes ((p^{1/2}\otimes 1)u^{c}(p^{-1/2}\otimes 1)))$ is a corepresentation of $\q$ , so we got 
%\begin{align*}
   % \Delta(u_{ij}u^{*}_{i_{1}j_{1}})=\sum_{k_{1},k_{2}} u_{ik_{1}}u^{*}_{i_{1}k_{2}}\otimes u_{k_{1}j}u^{*}_{k_{2}j_{1}}
%\end{align*}
%If we choose $a=e_{i_{0}j_{0}}$ then 
%\begin{align*}
  %  (Ad_{U}\otimes id)Ad_{U}(e_{i_{0}j_{0}})=&(Ad_{U}\otimes id)(\sum e_{ij}e_{i_{0}j_{0}}e_{j_{1}i_{1}}\otimes u_{ij}u^{*}_{i_{1}j_{1}})\\
    %=&(Ad_{U}\otimes id)(\sum _{i,i_{1}}e_{ii_{1}}\otimes u_{ii_{0}}u^{*}_{i_{1}j_{0}})\\
   %=&\sum_{i_{2}i_{3},i,i_{1}} e_{i_{2}i_{3}}\otimes u_{i_{2}i}u^{*}_{i_{3}i_{1}}\otimes u_{ii_{0}}u^{*}_{i_{1}j_{0}}\\
   % =&(Id\otimes \Delta)(\sum e_{i_{2}i_{3}}\otimes u_{i_{2}i_{0}}u^{*}_{i_{3}j_{0}})\\
   % =&(Id\otimes \Delta)Ad_{U}(e_{i_{0}j_{0}})
%\end{align*}
%Hence, $(Ad_{U}\otimes id)Ad_{U}(a)=(Id\otimes \Delta)Ad_{U}(a)$, for any $a\in B(H)$ and density conditions follows as  $(U\ot \bar{U_{r}})$ is a corepresentation of $\q$.
%\end{proof}

All the above results can be easily extended to left/ bi projective  corepresentation. In particular, we have
\blmma  Let $U\in B(H)\ot \q$ is a unitary. Then  followings are equivalent: 
     \begin{enumerate}
         \item [1)]If $U$ is a left projective  corepresentation.  \item[2)]There exists a positive matrix $\rho\in B(\bar{\clh})$ such that  $ ((\rho^{1/2}\otimes 1)U^{c}(\rho^{-1/2}\otimes 1))\ot U $  is a unitary corepresentation of $(\q,\Delta).$
         \item[3)]$U^{c}\ot U$ is an invertible corepresentation on $\bar{\clh}\ot\clh.$
     \end{enumerate} 
\elmma
\begin{proof} The proof of this is an easy consequence of the lemma \ref{cor}, hence ommited.\end{proof}
\blmma Let $U\in B(\clh)\ot \q$ be a unitary, then followings are equivalent:\begin{enumerate}
\item[1)]$U$ is a bi projective corepresentation of $\q$. 
\item[2)]There exists a positive matrix $\rho\in B(\bar{\clh})$ such that  $ ((\rho^{1/2}\otimes 1)U^{c}(\rho^{-1/2}\otimes 1))\ot U , U\otimes ((\rho^{1/2}\otimes 1)U^{c}(\rho^{-1/2}\otimes 1)) $   are unitary corepresentations of $(\q,\Delta).$ 
\item[3)] $U\ot U^{c},U^{c}\ot U$ are invertible corepresentations of $\q.$
\end{enumerate}
\elmma

 Proof of this lemma is straightforward and hence ommited.
\section{Projective envelope of a compact quantum group}\label{3}
We begin with defining a category.
  The construction of the category, the associated CQG or DQG, corepresentation $\clu$ and the homomorphism $\Pi$ in the proof of theorem \ref{main} are essentially taken from \cite{GS26}. The basic idea behind the construction of the category was motivated by Stefan Vaes through private communication to the first author of this paper.\\
 For a finite dimensional Hilbert space $\clh,$ we denote by $\bar{\clh}$ is the conjugate of $\clh.$ For $s\in \bar{\clh}\ot\clh, t\in\clh\ot\bar{\clh}$ ,  $R_s:\bbc\rightarrow \bar{\clh}\ot\clh, R_t:\bbc\rightarrow \clh\ot\bar{\clh}$ given by $R_s(1)=s$ and $R_t(1)=t.$\\
Let $\alpha$ be the category whose objects are $(U,\clh)$ where $\clh$ is a finite dimensional Hilbert space and the following conditions hold:
\begin {enumerate}
    \item[1)]$U\in B(\clh)\otimes \q$ is a unitary element.
    \item[2)]There exists a unitary element  $\bar U \in B(\bar{\clh} )\otimes \q$ and $s_U\in \bar{\clh}\ot\clh, t_{U}\in \clh\ot\bar{\clh}$  such that 
\begin{align}\label{1}
     &(R_{s_{U}}\otimes Id)(1_c\otimes Id)=(\bar {U}\otimes U)(R_{s_{U}}\otimes 1)  \\
     & (R_{t_{U}}\otimes Id)(1_c\otimes Id)=( U\otimes \bar {U})(R_{t_{U}}\otimes 1)
\end{align} and
\begin{align}\label{first equation}
  (R_{s_{U}}^* \otimes Id_{\bar {\clh}})(Id_{\bar {\clh}} \otimes R_{t_{U}})=Id_{\bar {\clh}}
 && (R_{t_{U}}^*\otimes Id_{\clh})(Id_{\clh}\otimes R_{s_{U}})=Id_{\clh}.
\end{align}
    \end{enumerate}

For any two objects $(U,\clh)$, $(V,\clk)$ of this category, the morphism space defined by
\begin{align*}
    &Mor((U,\clh),(V,\clk)))\\
    :=&\{T\in B(\clh,\clk):V(T\otimes 1)=(T\otimes 1)U~~~~,V\in B(\clk)\otimes \q,U\in B(\clh)\otimes \q\}.
\end{align*}
$U\ot V$ define by $U_{13}V_{23}.$ $U\oplus U$ in the usual way, as in \cite{GS26}.
\bthm
    $\alpha$ is a Rigid $ C^{\ast}$ tensor category.\label{22}
\ethm
\begin{proof}
This follows from Theorem 2.4 on \cite{LR97} as observed in \cite{GS26}.
\end{proof}

We have a canonical fiber functor
\begin{align*}
   & F:\alpha\rightarrow Hilb_{f}   ~~ {\mathrm {given~ by}}\\
   &F(U,\clh)=\clh ,F(T)=T,  \forall   T\in Mor((U,\clh),(V,\clk))\\
\end{align*}
By Tannaka-Krein duality, there is a CQG $Q_{univ}$ or equivalently the dual DQG $\hat{\q}_{univ}$ such that each $(U, \clh)\in obj(\alpha)$ is a unitary corepresentation of $\q_{univ}$. %From general theory of CQG corepresentation, it follows that there exists a positive matrix $\rho\in B(\bar{H})$ such that $(\rho\ot 1)u^{c}(\rho^{-1}\ot 1)$ is a unitary in $B(\bar{H})\ot \q.$
\blmma let $U\in B(\clh)\ot \q$ is a unitary such that there exists a positive invertible matrix $\rho\in B(\bar{\clh})$ for which  $(\rho\ot 1)U^{c}(\rho^{-1}\ot 1)$ is a unitary. Then $(U,\clh_{U})$ is an object of $\alpha.$
\elmma
%
%$(U,\clh)$ is an object of the category$\alpha$ if and only if there exists a positive invertible matrix $\rho\in B(\bar{H})$ such that  $(\rho\ot 1)u^{c}(\rho^{-1}\ot 1)$ is a unitary.\elmma
\begin{proof}
If $(U,\clh)$ is an object of this category $\alpha$ then $(\rho\ot 1)U^{c}(\rho^{-1}\ot 1)$ is a unitary, is already discussed above.
Let us assume that, $(\rho\ot 1)U^{c}(\rho^{-1}\ot 1)$ is a unitary in  $B(\bar{\clh})\ot \q.$ 
Now, defined the maps  $R_{s}:\bbc\rightarrow \bar{\clh}\ot \clh$ and $R_{t}:\bbc\rightarrow \clh\ot\bar{\clh}$ given by \begin{align*}
    &R_{s}(1)=\sum_{i}\rho^{-1}(\bar{e_{i}})\ot e_{i}\\
    &R_{t}(1)=\sum_{i} e_{i}\ot \rho(\bar{e_{i}}),
\end{align*}
where $\{e_{i}\}$ is an orthonormal basis of $\clh$ and $\{\bar{e_{i}}\}$ is the corresponding orthonormal basis for $\bar{\clh}$. Now, it is easy to prove it satisfies the following relations (\ref{1}), 2 and (\ref{first equation}). Hence $(U,\clh)$ is an object of this category $\alpha$.
\end{proof}

In our work, we will often deal with suitable $C^{*}$ subcategories of $\alpha$. To make it convenient for later use, let us give a general construction of such subcategories from a given family of unitary operators.
Let $C_{\alpha}$ be a subset of $Obj(\alpha)$ and the objects of $C_{\alpha}$ are of the form $(U,\clh).$ Let $\bar{C}_{\alpha}$ the smallest subfamily of $Obj(\alpha)$ containing $C_{\alpha}$ which is closed under direct sums and tensor products and subobjects.

%Suppose , we choose $C_{\alpha}$ is the set of all biprojective corepresentation. Then $\tilde{C_{alpha}}=Biproj(C)$ is  the category generated by all irreducible biprojective corepresentation.
%objects of $Biproj(C)$ are in the form $((p\ot 1)u,H)$ and $p\in End(u)$ is a projection and $u$ is an object of $Biproj(C)$ and ,\begin{align*}
  %  &Mor(((p\ot 1)U,H),((q\ot 1)V,K))=qMor(U,V)p,  p\in End(u),q\in End(V)\\
    %&\mathrm{p,q~boths~are~projections}.
%\end{align*} So, $Biproj(\tilde{C})$ is a strict monoidal category and it is closed under direct sum and  subobject from previous theorem (\ref{22}).\\

%Let $F_{Biproj(C)}:Biproj(C)\rightarrow Hilb_{f}$ be a fiber functor, is given by the restriction of the fiber functor on the category $Biproj(C)$.
   
%\begin{align*}
   % &F(u,p)=pH_{u} ,u\in B(H_{u})\otimes Q\\
     %&F:Mor((u,p),(V,q))\rightarrow Mor(pH_{u},qH_{v}) \\
    % &s.t. F(T)=T,T\in Mor((u,p),(V,q))
%\end{align*}
%$F$ is linear on morphisoms and $F_{2}$ is the natural identity map and $F(T^{*})=F(T)^{*}=T^{*}$. so, $F$ is a unitary exact fiber functor.\\
%\textbf{Woronowicz's Tannaka-Krein duality:}\\
%Let $C$ be a $C^{*}$ tensor category with conjugates.$F:C\rightarrow Hilb_{f}$ be a unitary fiber functor. Then there exist a Compact quantum group $Q$ and unitary monoidal equivalence $E:C\rightarrow Rep(Q)$ s.t. $F$ is naturally unitary monidally isomorphic to the composition of the canonical fiber functor $Rep(Q)  \rightarrow Hilb_{f}$ with $E$.\\
 Let us consider the following choices of $C_{\alpha}$:\begin{enumerate}
     \item [1)]$C^{L}_{\alpha}= \{(U,\clh):\mathrm{U~ is ~a~left~projective~corepresentation~of ~\q~on~\clh}\}.$
     \item[2)]$C^{R}_{\alpha}= \{(U,\clh):\mathrm{U~ is ~a~right~projective~corepresentation~of ~\q~on~\clh}\}.$
     \item[3)]$C^{bi}_{\alpha}= \{(U,\clh):\mathrm{U~ is ~a~bi~projective~corepresentation~of ~\q~on~\clh}\}.$
 \end{enumerate}

     Apply the above construction, we will get the categories $\bar{C}^{L}_{\alpha}\equiv\bar{C}^{R}_{\alpha},\bar{C}^{bi}_{\alpha}.$
% So, we can say that there exist a compact quantum group $\tilde{Q}$ and a unitary monoidal equivalence $E:Bipro(\tilde{C})\rightarrow Rep(\tilde{Q})$ s.t. $F$ is naturally unitary monoidal isomorphic to the composition of the canonical fiber functor $Rep(Q)\rightarrow Hilb_{f}$ with $E$.\\
 % Let, $Q=(\q,\Delta)$ be a CQG. Now, we take the $C^{*}$ tensor category $Rep(Q)$, where objects are finite dimensional unitary corepresentation and morphisoms between two objects are just the intertwiners of those two corepresentations. We know every finite dimensional corepresentation is a projective corepresentation of that compact quantum group $Q$.
 % Now, we define a functor \begin{align*}
   %  & F_{1}:Rep(Q)\rightarrow Biproj(\tilde{C})\\
   %  & s.t. F_{1}(u)=(u,Id_{H_{u}})\\
    %  &F_{1}:Mor(u,V)\rightarrow Mor((u,id_{u}),(V,Id_{v}))\\
      %&F_{1}(T)=T
 % \end{align*}\\
  %$F_{1}$ is a fully faithful unitary tensor functor or we can say that $Rep(Q)$ is an embeded subcategory inside biprojective $C^{*}$ tensor category of $Q$.
   \bthm
   \label{main}
       Let $(\q,\Delta)$ be a CQG. Then there exist CQGs $(\widetilde{\q_{L}},\Delta_L)\equiv(\widetilde{\q_{R}},\Delta_R),\\ 
 (\tilde{\q_{bi}},\Delta_{bi})$ along with injective CQG morphisms \begin{enumerate}
           \item [1)]$i_{L}:\q\to \widetilde{\q_{L}},$
           \item[2)]$i_{R}:\q\to \widetilde{\q_{R}},$
           \item[3)] $i_{bi}:\q\to \widetilde{\q_{bi}}
        $ 
       \end{enumerate} such that identifying $\q$ as a Woronowicz's subalgebra of each of these CQGs $\widetilde\q_{L}\equiv\widetilde\q_{R},\widetilde\q_{bi},$
            the following properties hold:\begin{enumerate}
           \item [i)] For every left projective unitary corperesentation $(U,\clh)$ of $\q$ there exists a unitary (linear) corepresentation $(\tilde{U}_{L},\clh)$ of $\widetilde{\q}_{L}$ such that $Ad_{\tilde{U}_{L}}$ is isomorphic to $Ad_{U}.$ A similar statement holds for right/bi projective unitary corepresentation of $\q.$
           \item[ii)] The *-algebra generated by the matrix coefficients of $\tilde{U}_{L}$ is same as the Hopf $^{*}$-algebra $Pol(\widetilde{\q}_{L})$, where $U$ varies over the set of all left projective unitary corepresentations of $\q.$ A similar statement holds for right/bi projective unitary corepresentation of $\q.$
       \end{enumerate}
   
   \ethm
   \begin{proof}
We prove the statements (i) and (ii) only for $\bar{C}^{L}_{\alpha}$, as the proof the other case are very similar. %The UTC $\widetilde{C^{L}_{\alpha}}$ are the Tannka-Krein construction 
Let $\q\subseteq B(\clh_0),$ where $\clh_0$ is the GNS space for the Haar state on $\q.$ Recall the Tannaka-Krein construction, which involves the $*$-algebra $Nat(F^{L})$, where $F^{L}$ is the restriction of the canonical fiber functor $F$ to the subcategory $\bar{C^{L}_{\alpha}}.$ Any element $T$ of $Nat(F^{L})$ is given by a collection $\{T_{(U,\clh_{U})}: T_{(U,\clh_{U})}\in B(\clh_{U})\}.$ 

The von-Neumann algebraic DQG  $M$ of $\widetilde{\q_L}$ is isomorphic with 
\begin{align}
    M=\{T\in Nat(F^{L}):sup _{(U,H)\in \text{Obj}(\bar{C}^{L}_{\alpha}})\|T_{(U,H_U)}\|<\infty\}
\end{align}
For $\xi,\eta\in \mathcal{H}_0$, consider $\mathcal{U}_{\xi,\eta}\in \text{Nat}(F^L)$ given by 
$${\big(\mathcal{U}_{\xi,\eta}\big)}_{(U,H_U)}=U_{\xi,\eta}$$ where $U_{\xi,\eta}\in \mathcal{Q}$ is the contraction of $U\in \mathcal{B}(\mathcal{H}_U)\otimes \mathcal{Q}$, i.e. $$\langle \theta_1,U_{\xi,\eta}(\theta_2)\rangle=\langle \theta_1\otimes \xi, U(\eta\otimes \theta_2)\rangle \text{ for all } \theta_1,\theta_2\in \mathcal{H}_0. $$
As each $U$ is unitary, observe that $\text{sup}_{(U,\mathcal{H}_U)\in \text{Obj}(\widetilde{C}^L_{\alpha})}\|U_{\xi,\eta}\|< \|\xi\|\|\eta\|$, which proves that $\mathcal{U}_{\xi,\eta}\in \mathcal{M}.$\\
Let $\{(U_i,\mathcal{H}_i):i\in I\}$ be an enumeration of all the irreducible objects of the category $\bar{C}^L_{\alpha}$.
Then $\mathcal{M}\cong\prod_{i\in I}\mathcal{B}(\mathcal{H}_{U_i})\text{ and } \mathcal{M}\subseteq \mathcal{B}(\oplus_{i\in I}\mathcal{H}_{U_i}),$  so that $\mathcal{M}\otimes \mathcal{B}(\mathcal{H}_0)\subseteq \mathcal{B}(\oplus_{i\in I}(\mathcal{H}_{U_i}\otimes \mathcal{H}_0)$.
\\Viewing $\mathcal{U}$ in this picture we see that $\mathcal{U}=\prod_{i\in I}U_i,$ where $U_i\in \mathcal{B}(\mathcal{H}_{U_i})\otimes \mathcal{Q}\subseteq \mathcal{B}(\mathcal{H}_{U_i}\otimes \mathcal{H}_0).$
As each $U_i$ is unitary, so is $\mathcal{U}$. In Lemma 3.5 of \cite{GS26}, we define  $\Pi:=\Pi_{\clu}$ from the Hopf *-algebra $(\widetilde{\q^{L}})_0$, spanned by the matrix elements of irreducible corepresentations of $\widetilde{\q_{L}},$ to $\q\subseteq B(\clh_0)$ given by
\begin{align}
    \Pi(\omega)=(\omega \ot Id)(\clu),
\end{align} where $\omega\in (\widetilde{\q_{L}})_0.$ This is well defined by the general theory of CQG/DQG and $\Pi$ is a $\ast$-homomorphism. For $i\in I,$ let $\widetilde{U_{i,L}}$ be the irreducible unitary corepresentation of $\widetilde{\q_{L}}$ on $\clh_{U_i}.$ It is easy to see that $(Id\ot \Pi)(\widetilde{U_{i,L}})=U_{i}.$

In fact, in a similar way we can show that $(Id\ot \Pi)(\widetilde{U_{L}})=U$, for all objects $(U,\clh_{U})$, where $\widetilde{U_{L}}$ is the corresponding unitary corepresentation of $\widetilde{\q_{L}}$ on $\clh_{U}.$\\

Now, note that if $(W,\clh_w)\in obj(Corep(\q))$ that is $W$ is actually a unitary corepresentation of $\q$ on $\clh_{W},$ then $\Pi(\omega)=\omega$ for any matrix coefficient $\omega$ of $W.$ In other words $(Id\ot \Pi)(W)=W.$ One can also obeserve that
\begin{align*}
    \Pi(\widetilde{U^c_{L}})=(\Pi(\widetilde{U_{L}}))^c=U^c
\end{align*}, for all objects $(U,\clh_{U}).$

By construction, for all objects $(U,\clh_U)$, $(U^c\ot U,\bar{\clh}\ot \clh)$ is a corepresentation of $\q.$
But the unitary corepresentation to the object $(U^c\ot U,\bar{\clh}\ot \clh)$ is nothing but $\widetilde{U^c_{L}}\ot \widetilde{U_{L}}.$ Hence by our earlier observation 
\begin{align}
   (Id\ot \Pi)(\widetilde{U^c_{L}}\ot \widetilde{U_{L}})=(\widetilde{U^c_{L}}\ot \widetilde{U_{L}})
\end{align} so $U^c\ot U=\widetilde{U^c_{L}}\ot \widetilde{U_{L}}$ which implies that $\widetilde{U_{L}}$ is a left projective  corepresentation on $\clh_{U}.$ This proves (i).\\

(ii) Follows from the construction and the general theory of CQG.
\end{proof}
\bcrlre let $U\in B(\clh)\ot \q$ is a unitary. $(U,\clh)$ is an object of the category $\alpha$ if and only if there exists a positive invertible matrix $\rho\in B(\bar{\clh})$ such that  $(\rho\ot 1)U^{c}(\rho^{-1}\ot 1)$ is a unitary.
\ecrlre
We call the quantum groups $\widetilde{\q_{L}},\widetilde{\q_{R}}, \widetilde{\q_{bi}}$ are left, right, bi projective envelope of $\q.$
\section{Normal quantum subgroups and normalizer of a quantum subgroup}
As before $(\q,\Delta)$ be a CQG with the (mutually inequivalent) irreducible uitary corepresentation $U^{\alpha}=\sum^{d_{\alpha}}_{i,j=1} u^{\alpha}_{ij}\otimes e^{\alpha}_{ij}$, where $\alpha\in L$ ($L$ is an indexed set). Let $H\subseteq L$ and $\q_{H}$ be a Woronowicz's subalgebra of $\q$ and  $\hat{\q}_{H}$ be the dual DQG of $\q_{H}.$

We recall the notion of normal quantum subgroups of a DQG and quotient with respect to them from \cite{normalvai}, \cite{normalalex} . %We need the special case of normal quantum subgroups of DQG, which we discuss in some details now% .
    Let $\hat{\clq}$ be a DQG given by $C_{0}$ direct sum of $\clb(\clh_\alpha), \alpha \in I$, where $d_\alpha={\rm dim}(\clh_\alpha)$ and let $\clq$ be the dual CQG with the mutually inequivalent  irreducible corepresentations $u^\alpha $ corresponding to $\alpha  \in I$, given by $(( u^\alpha_{ij} )) \in \clb(\clh_\alpha) \ot \clq$. Let $W=\sum_{\alpha , i,j} (u^\alpha_{ji})^* \ot e^\alpha_{ij} 
     \in (\hat{\clq} \ot \clq)''$  be the so-called left regular corepresentation, where $e^\alpha_{ij}, i,j=1, \ldots, d_\alpha$ are the `matrix units' of $\clb(\clh_\alpha)$.

     A quantum subgroup of $\hat{\clq}$ is given by a both sided weakly closed Hopf ideal. Any such ideals of a $C_{0}$ direct sum of matrix algebras is again one such a $C_{0}$ direct sum, say $\oplus_{ \alpha \in J} \clb(\clh_\alpha)$ for some subset $J $ of $I$. Then the quotient is again a $C_{0}$ direct sum, over $H:=J^{(c)}$ (complement of $J$). We write $\hat{\clq}_J$ and $\hat{\clq}_H$ for the $C_{0}$ direct sum over the index sets $J$ and $H$ respectively and denote by $P_J$, $P_H$ the corresponding central projections in $\hat{\clq}$ onto $\hat{\clq}_J$, $\hat{\clq}_H $ respectively. The co-ideal property of $\hat{\clq}_J$ implies $(P_H \ot P_H) \Delta(a)=0$ for all $a \in \hat{\clq}_J$. 
     This implies that for $\alpha, \beta \in H$ and for any $\gamma \in I$ such that the irreducible index by $\gamma$ is a direct summand of the tensor product of the irreducible corepreentations corresponding to $\alpha$ and $\beta$, we must also have $\gamma \in H$ too. 
     
    % This means $\{ u^\alpha_{ij}, \alpha \in H,~i,j=1, \ldots, d_\alpha\}$ forms a unital $\ast$ subalgebra of $\hat{\clq}$. We denote the weak closure of this algebra in $\hat{\clq}$ 
 %by $\hat{\clq_H} $.
 
Note that $\hat{\clq}_H $ has the coproduct given by $\Delta_H:=(P_H \ot P_H) \circ \Delta$ restricted to $\hat{\clq}_H$. Indeed, it is clear that $\Delta_H$ maps $\hat{\clq}_H$ to $\hat{\clq}_H \ot \hat{\clq}_H$ and the coassociativity of the map follows from the coassociativity of $\Delta$ combined with the fact that $\Delta_H((1-P_H)(a))=0 $, implying $\Delta_H(P_Ha)=\Delta_H(a) $ for all $a$. There is a surjective quantum group morphism from $\hat\clq$ to $\hat\clq_H$ which sends $a$ to $P_H(a)$. 
 \bdfn\label{norm}

$\hat {\q}_{H}$ is said to be normal quantum subgroup of $\widehat{\q}$ if $W_{\q}(\q_{H} \otimes Id)W^{*}_{\q} \subset \q''_{H}\bar{\otimes} B(L^{2}(\q)),$ where $L^2(\q)$ denotes the  GNS space of $\q$ with respect to the Haar state and $W_{\q}$ denotes the left regular representation of $\q,$ which is given by $W_{\q}=\sum_{\alpha,i,j} u^{\alpha}_{ij} \ot e^{\alpha}_{ij}.$
\edfn
 % We take the characterization of normality of a quantum subgroup as in \cite{vaes_normal} and \cite{alex_normal} as the definition given below:
%\bdfn
%\label{normaldef}
% The quantum subgroup $\hat{\q}_H$ is said to be normal if $W(\clq_H  \ot 1)W^* \subseteq \clq_H  \ot \clb(\clh)$, or equivalently, $W (\clq_H \ot 1) \subseteq (\clq_H \ot \clb(\clh)) W$. Here, $\clh$ denotes the Hilbert space %$\bigoplus_{\alpha \in I} \clh_\alpha \ot \overline{\clh_\alpha}$, 
% which is isomorphic with the GNS space of $\clq$  w.r.t. its Haar state. 
 %\edfn
 It is well-known (see \cite{normalalex}, \cite{normalvai}) and references therein that 
 \bppsn 
 $\hat\clq_H$ is said to be normal subgroup of $\hat\clq$ if and only if the left-invariant subalgebra $\{ a \in \hat{\clq} :~ ((P_H \ot {\rm id} )\circ \Delta)(a)=1_H \ot a \}$ coincides with the right-invariant subalgebra $\{ a \in \clq :~ (({\rm id} \ot P_H) \circ \Delta) (a)=a \ot 1_H  \}$.
 \eppsn
 \bcrlre If $\hat\clq_{H}$ is a normal subgroup of $\hat\clq$ then the left/right invariant subalgebra is a subset of $ker(P_{H})\cup \{\bbc 1_{\clq_{H}}\}$ .\ecrlre 
% \begin{proof}
 %Let $a\in\hat{\clq_{H}}$ such that $((P_H \ot {\rm id} )\circ \Delta)(a)=1_H \ot a \}.$ It will also satisfy $(({\rm id} \ot P_H) \circ \Delta) (a)=a \ot 1_H $ and \begin{align*}
  %   (P_{H}\ot P_{H})\Delta(a)=& (P_{H}\ot Id)(Id\ot P_{H})\Delta(a)\\
   %  =&(P_{H}\ot Id)(a\ot 1_{H})\\
    % =& P_{H}(a)\ot 1_{H}.
 %\end{align*}
%Also, \begin{align*} 
 %  (P_{H}\ot P_{H})\Delta(a)=&(Id\ot P_{H})(P_{H}\ot Id)\Delta(a)\\
  % =&(Id\ot P_{H})(1_{H}\ot a)\\
  % =&1_{H}\ot P_{H}(a)
 %\end{align*} So, $P_{H}(a)=0$ or $P_{H}(a)=c1_{H}.$
 
 \vspace{2mm}  
 
   If $\hat{\clq}_H$ is a normal subgroup, the corresponding %(von Neumann algebraic) 
   quotient quantum group, say $\hat{\clq}/\hat{\clq}_H$, is given by the left (equivalently right) invariant subalgebra mentione above,  which also coincides with   the  $C^{\ast}$-algebra generated by the image of the operator-valued weight $T_H$ given by $T_{H}(\cdot)=(\tau_H \ot {\rm id}) \circ \Delta$ on $\hat\clq \subseteq \clb(\clh)$ where 
 $\tau_H$ denotes the restriction of $\tau$ to $\hat\clq_H$, which is a   semifinite operator valued weight with all $\clq_\alpha$'s in the domain in particular.

 \blmma
 \label{normallemma}
 Suppose for each $\alpha \in I$ and $\beta \in H$, there is some  $V^{\alpha,\beta} \in \clq_H \ot \clb(\clh_\alpha \ot \clh_\beta)$ such that $(u^\alpha)_{(12)} (u^\beta)_{(13)} =V^{\alpha, \beta}(u^\alpha)_{(12)}$. Then $\hat{\clq}_H$ is a normal quantum subgroup of $\clq$.
 
 \elmma
 \begin{proof}

 It suffices to verify the conditions of Definition \ref{norm} for the generators $u^\beta_{kl} $ of $\clq_H $, for $\beta \in H$. But from assumption, we have $W_{(12)} (u^\beta)_{(13)}=VW_{(12)}$, with $V:=\bigoplus_{\alpha \in I} V^{\alpha, \beta} \in \hat{\clq}_H \ot \clb(\clh_\alpha) \ot \clb(\clh_\beta).$ Taking $\phi$ to be the functional on $\clb(\clh_\beta)$ which is $1$ on $e^\beta_{kl}$ and $0$ on other matrix elements of $\clb(\clh_\beta)$ and applying $({\rm id }_{(12)} \ot \phi)$ on both sides of the above equation, the proof of the lemma is complete.  \end{proof}

 \blmma\label{4.2} $\hat{\q}_H$  is normal if and only if for any $x\in {\q_{H}}$  and $\alpha\in L, 1\leq i,i_1\leq d_{\alpha},$ $v^{\alpha,x}_{ii_{1}}=\sum_{j}u^{\alpha}_{ij}xu^{*\alpha}_{i_{1}j}$ is in ${\q_{H}}$\elmma
 \begin{proof}
 We note that for any $x$ in $\q_{H}$,
 \begin{align*}
     W_{\q}(x\otimes 1)W^{*}_{\q}=&(\sum_{\alpha}(\sum_{i,j}u^{\alpha}_{ij}\otimes e^{\alpha}_{ij}))(x\otimes 1)(\sum_{\alpha}(\sum_{i_{1},j_{1}}u^{*\alpha}_{i_{1}j_{1}} \otimes e_{i_{1}j_{1}}))\\
     =&\sum_{\alpha}(\sum_{i_{1},j_{1},i,j}u^{\alpha}_{ij}xu^{*\alpha}_{i_{1}j_{1}}\otimes e^{\alpha}_{ij}e^{\alpha}_{j_{1}i_{1}}\\
     =&\sum_{\alpha}(\sum_{i,i_{1}}u^{\alpha}_{ij}xu^{*\alpha}_{i_{1}j}\otimes e^{\alpha}_{ii_{1}})\\
     =&\sum_{\alpha}\sum_{i,i_{1}} v^{\alpha,x}_{ii_{1}} \otimes e^{\alpha}_{ii_{1}}.
 \end{align*}
Let us assume that $\hat{\q}_{H} $ be a normal quantum subgroup of $\widehat{\q} $. Then $\sum_{i,i_{1}} v^{\alpha,x}_{ii_{1}} \otimes e^{\alpha}_{ii_{1}}$ is in $\q''_{H}\bar\otimes B(L^{2}(\q))$, if we take the normal linear functional $\phi^{\alpha}$ on $B(L^{2}(\q))$ such that $ \phi^{\alpha}(e^{\alpha}_{i_{0}j_{0}})$=1 and $\phi^{\alpha}(e^{\alpha}_{i_{0}j_{0}})$=0 if $(i ,j)\neq(i_{0}, j_{0})$ . Then (Id$\otimes$ $\phi^{\alpha}$) ($W_{\q}(x\otimes 1)W^{*}_{\q}$)=$v^{\alpha,x}_{i_{0}j_{0}}  \in \q_{H}.$ Now, for $x\in Pol(\q_{H}),$
clearly $V^{\alpha,x}_{i_0j_0}\in Pol(\q).$ Hence it is in $Pol(\q)\cap \q''_{H}=Pol(\q_{H}).$ In general approximating  $x\in \q_{H}$ by a sequence from $Pol(\q_{H}),$ we get $v^{\alpha,x}_{ii_{1}}\in \q_{H}.$

The converse also follows from the Lemma \ref{4.2} as the sum is a strongly convergent series and $\q''_{H}\bar\ot B(L^{2}(\q))$ is a von-Neumann algebra. 
\end{proof}
\brmrk\label{4.3}
It is clear from the proof of Lemma \ref{4.2} that it is enough to verify the condition of the lemma for any subset $K\subseteq L$ which is generating in the sense that every object in $L$ is a subobject of finite direct sums of tensor products of objects from the subset $K.$
\ermrk
Now prove a result in the categorical setting, which will be useful for our later discussion.
Let us assume that $C$ is a UTC and $F$ is a fiber functor on $C$. Suppose also that $C_{0}$ is an embedded sub UTC of $C$ and $F_{0}$ is the restriction  of $F$ on $C_{0}$. Let $\{x_{\alpha}:\alpha\in I\} $ be collection of irreducible objects of $C$ and $\{x_{\alpha}:\alpha\in I_{0}\} $ be the collection of irreducible objects of $C_{0}.$\\
Let $A_{J}$=\{$x\in obj(C):x \otimes y\otimes \bar{x},\bar{x}\otimes y\otimes x\in obj(C_{0}),\forall y\in obj(C_{0})$ \}.\\
Now, we define a subcategory $C_{J}$ by $obj(C_{J})=\{\bigoplus^{n}_{i=1}{x_{i}}: x_{i}\in A_{J},n\geq 1\}.$ %where
%\begin{align*}
   % x_{i}\otimes y\otimes \bar{x_{i}}\in obj(C_{0}) ,           ~~~\bar{x_{i}}\otimes y\otimes x_{i}\in obj(C_{0}).
%\end{align*}
\par
\blmma $C_{J}$ is a UTC and $C_{0}$ is a sub UTC of $C_{J}$.\elmma
\begin{proof}
Let $\bigoplus{x_{i}},\bigoplus{y_{j}}\in obj(C_{J})$ (finite direct sums). Now ($\bigoplus{x_{i}})\otimes (\bigoplus{y_{j}}) $  is isomorphic to $\bigoplus(x_{i}\otimes y_{j}).$ For any object $z$ of $C_{0}$, $x_{i}\otimes y_{j}\otimes z\otimes\bar{y_{j}}\otimes\bar{x_{i}}$ is an object of $C_{0}$ because  $y_{j}\otimes z\otimes\bar{y_{j}}\in obj({C_{0}})$ and hence $x_{i}\otimes (y_{j}\otimes z\otimes\bar{y}_{j})\otimes\bar{x}_{i}\in Obj(C_0)$ too. This means $x_i\otimes y_i\in A_{J},$ hence their finite sum is also an object of $C_{J}.$
%So, $x_{i}\otimes (y_{j}\otimes z\otimes\bar{y_{j}})\otimes\bar{x_{i}}$ is an object of $C_{0}$. Direct sum of objects is an object of $C_{J}$ that follows trivially from construction of $C_{J}$.\\

Next we will prove it is closed under taking subobject. Let $x_{p}$ be a subobject of $x$, where $x\in obj(C_{J}).$ We have $x \otimes y\otimes \bar{x},\bar{x}\otimes y\otimes x\in obj(C_{0}).$ So, there exist an isometry $V_{p}\in Mor(x_{p},x)$ and $V_{p}V^{*}_{p}=P\in End(x)$. We know if $x_{p}$ is a subobject of $x$ then $\bar{x}_{p}$ is a subobject of $\bar{x}$. Similarly, we can say that there exist an isometry $\bar{V_{p}}\in Mor(\bar{x}_{p},\bar{x})$ and $\bar{V_{p}}\bar{V^{*}_{p}}=\bar{P}\in End(x)$. Hence it follows that $(V_{p}\otimes Id_{y}\otimes \bar{V}_{p})\in Mor(x_{p}\otimes y\otimes \bar{x}_{p},x\otimes y\otimes \bar{x})$ and it is an isometry and $V_{p}V^{*}_{p}\otimes Id\otimes \bar{V_{p}}\bar{V^{*}_{p}}=(p\otimes Id\otimes \bar{p})\in End(x\otimes y\otimes \bar{x})$. So, $x_{p}\otimes y\otimes \bar{x_{p}}$ is a subobject of $(x\otimes y\otimes \bar{x})\in obj(C_{0})$ and as $C_{0}$ is closed under taking subobjects,  $x_{p}\otimes y\otimes \bar{x}_{p}\in obj(C_{0}).$  Now, if we choose an object $x=\bigoplus^{n}_{i=1} x_{i}$ then any subobject is of the form $\bigoplus x_{p_{i}}$  ,where $x_{p_{i}}$ is a subobject of $x_{i}$ as we know $x_{p_{i}} \otimes y\otimes \bar{x}_{p_{i}},\bar{x}_{p_{i}}\otimes y\otimes x_{p_{i}}\in obj(C_{0}).$ This implies $x\in obj(C_{0}).$ It follows easily that $C_{J}$ is also closed under taking conjugate. Finally, we note that $obj(C_{0})\subseteq obj(C)$ from the definition of $C_{J}$ and the fact that $C_{0}$ is an embeded sub UTC.
%and rigidity of this category follows easily. So,$C_{j}$ is a uTC.\\
%It's clear that for any two objects $x,y$ of $C_{0}$, $x\otimes y\otimes \bar{x}$ is an object of $C_{0}$ and morphisom between two objects $x,y$ of $C_{0}$ is same as the $Mor_{C}(x,y)$. So,$C_{0}$ is a sub uTC of $C_{J}.$
\end{proof}

\bdfn\label{nz} We call $C_{J}$ the normalizer of $C_0$ in $C$ and define it by $N(C,C_0).$
\edfn
 Applying Tannaka-krein duality on $C_{J}$ and $C_{0}$ with the fiber functor $F$ and $F_{0}$ respectively, we get CQG $\q_{J}$ and $\q_{0}$ such that $\q_{0}$ is a Woronowicz subalgebra of $\q_{J}.$ \\
 
The following theorem justify the terminology of definition \ref{nz}.

\bthm$\hat{\q}_{0}$ is a quantum normal subgroup of $\hat{\q}_{J}$.\ethm\label{normal}
\begin{proof}
We apply Lemma \ref{4.2}, as before, for an irreducible object $z$ in $A_{J},$ choose and fix a unitary corepresentation $U^{z}.$ By Remark \ref{4.3} it is enough to verify  the condition of lemma \ref{4.2} for the generating subset $A_{J}.$ If $x_{\alpha}$ is an irreducible object in $A_{J}$ and $x_{\beta}$ is an irreducible object of $C_{0}$ then $x_{\alpha}\otimes x_{\beta}\otimes \bar{x}_{\alpha}$ is an object of $C_{0}.$ So, there exists a unitary operator $T$ such that $(1_{\q}\otimes T)(U^{x_{\alpha}}\otimes U^{x_{\beta}}\otimes U^{x_{\alpha}})(1_{\q}\otimes T^{*})$ is a unitary coreprsentation of $\q_{0}$. Write $U^{x_{\alpha}}=\sum_{i,j} u^{\alpha}_{ij}\otimes e^{\alpha}_{ij}$, $U^{x_{\beta}}=\sum_{k,l} v^{\beta}_{kl}\otimes f^{\beta}_{kl}$,   $U^{\bar{x}_{\alpha}}=\sum_{m,n} u^{*\alpha}_{mn}\otimes R^{1/2}\bar{e^{\alpha}_{nm}}R^{-1/2}$ for suitable matrix units $e^{\alpha}_{ij},f^{\beta}_{kl} $ and positive invertible operator $R\in B(\bar{\clh_{\alpha}}).$\\
Let us assume that $(U^{x_{\alpha}}\otimes U^{x_{\beta}}\otimes U^{\bar{x}_{\alpha}})=V^{\alpha,\beta,\bar{\alpha}}=\sum u^{\alpha}_{ij} v^{\beta}_{kl} u^{*\alpha}_{mn}\otimes  e^{\alpha}_{ij}\otimes f^{\beta}_{kl}\otimes R^{1/2}\bar{e^{\alpha}_{nm}}R^{-1/2}$. Now, we take the linear functional $\phi^{\alpha}$ on $B(\clh_{\alpha})$ satisfy $\phi^{\alpha}(e^{\alpha}_{i_{0}j_{0}})$=1 and $\phi^{\alpha}(e^{\alpha}_{ij})=0$ if $(i,j)\neq (i_0, j_0)$. Similarly, we choose the linear functionals $\phi^{\beta}$ and $\phi^{\bar{\alpha}}$ on $B(\clh_{\beta})$,$B(\bar{\clh_{\alpha}})$ given by  $\phi^{\beta}(e^{\beta}_{k_{0}l_{0}})=1$,~$\phi^{\beta}(e^{\beta}_{kl})=0$ if $(k, l)\neq (k_0,l_0)$ and  $\phi^{\bar{\alpha}}(R^{1/2}\bar{e^{\alpha}_{n_{0}m_{0}}}R^{-1/2})=1$ and   $\phi^{\bar{\alpha}}(R^{1/2}\bar{e^{\alpha}_{nm}}R^{-1/2})=0$ if $(n,m)\neq (n_0,m_0)$.\\
    Then,$(Id\otimes \phi^{\alpha}\otimes \phi^{\beta}\otimes \phi^{\bar{\alpha}})(V^{\alpha,\beta,\bar{\alpha}})=u^{\alpha}_{i_{0}j_{0}}v^{\beta}_{k_{0}l_{0}}v^{*\alpha}_{m_{0}n_{0}}$ % So, for any fixed $v^{\beta}_{k_{0}l_{0}}$, 
 %$\sum_{j}u^{\alpha}_{ij}v^{\beta}_{k_{0}l_{0}}u^{*\alpha}_{i_{1}j}$ 
 is in $Pol({\q_{0}})$. This proves that $\sum_{j} u^{\alpha}_{ij}v^{\beta}_{k_{0}l_{0}}v^{*\alpha}_{i_1j}$ is in $ Pol(\q_{0}),$%As any $x\in Pol(\q_0)$ is a limit of  finite sums of elements of the form $v_{k_0l_0}$, we can conclude that $\sum_{j} u^{\alpha}_{ij}xv^{*\alpha}_{i_1j}\in \q_{0}.$
 Finally, approximating a general element $x\in\q_{0}$ by a sequence from  $Pol(\q_{0}).$  we complete the proof.\\

 %Let us now recall the framework of unitary tensor category (UTC for short) or $C^*$ tensor category, from   which we refer the reader to \cite{snbook} and the references therein.
 %For two objects $x,y$ of some UTC, let us write $x \leq y$ to mean that $x$ is isomorphic with a subobject of $y$. 
  Consider a normal quantum subgroup $\hat{\clq}_H$ of $\hat{\clq}$ as in the preceeding discussion, with $\clc=\hat{\clq}/\hat{\clq}_H$ being the quotient DQG, identified as a Woronowicz subalgebra of $\hat{\clq}$. For a DQG $\cld,$ let ${\rm Rep}(\cld)$ denote the UTC of finite dimensional $\ast$-representations of $\cld$, which is also equivalent to the category ${\rm Corep}(\hat{\cld})$ of finite dimensional corepresentations of the dual CQG $\hat{\cld}$. Recall the set $I$ that labels the irreducibles of ${\rm Corep}(\clq)$.  Let $e$ denote the unit object, which corresponds to the counit $\epsilon$ of the DQG $\hat{\clq}$. The counit of $\clc$ is the restriction of $e$ to $\clc$ and  we denote the unit object 
   of ${\rm Rep}(\clc)$ by $e_\clc$.  We define a relation %(see \cite{} )  
   $\alpha \sim \beta $ for $\alpha, \beta \in {\rm Rep}(\hat{\clq})$  by:
   \bdfn\label{er}
   $\alpha \sim \beta$ if there is some  $h \in H$ such that ${\rm Mor}(\alpha,  \beta \ot h) \ne \{ 0 \}$ 
   \edfn

   %\bdfn
  % $\alpha \sim \beta$ if there is some  $h \in H$ such that ${\rm Mor}(\alpha,  \beta \ot h) \ne \{ 0 \}$ 
  % \edfn
  
  \blmma\label{4.11}
  
  (i) $\sim$ is an equivalence relation. \\
  (ii) for $\alpha, \beta \in I={\rm Irr}(\hat{\clq})$, $\alpha \sim \beta$ if and only if there is some $h \in H$ such that $\alpha \leq \beta \ot h$.
  \elmma
  The above lemma  follow from  the 4 th section of \cite{kenbs}.
  %(i) 
% If   ${\rm Mor}(\alpha, \beta \ot h) \ne \{ 0 \} $ we have  ${\rm Mor}(\alpha \ot \overline{h} , \beta)  \ne \{ 0 \}$, hence $\beta \sim \alpha$ which proves symmetry of the relation.  Other properties (reflexivity and transitivity) are also easy to prove, using the fact that $h \ot h^\prime $ is a direct sum of irreducibles in $H$ for $h,h^\prime \in H$.\\
 %(ii) This follows from the fact that for an irreducible $\alpha$ and any $x \in {\rm Rep}(\hat{\clq})$,  ${\rm Mor}(\alpha, x)$ is nonzero if and only if $\alpha $ is a sub-object of $x$, i.e. $\alpha \leq x$.\end{proof}
  
  There are canonical functors $F_1 : {\rm Rep}(\hat{\clq}_H) \raro {\rm Rep}(\hat{\clq})$ and $F_2 : {\rm Rep}(\hat{\clq}) \raro {\rm Rep}(\clc)$ given by the following:\\
  $$ F_1((\rho,\clh) )=(\rho \circ P_H, \clh), ~~~F_2((\theta, \clk))=(\theta|_{\clc}, \clk),$$
   where $P_H :\hat{\clq}  \raro \hat{\clq}_H$ the surjective map  and $\theta|_\clc$ denotes the restriction of $\theta$ to the subalgebra $\clc$ of $\hat{\clq}$. Both $F_1$ and $F_2$ are identity on morphisms. It is easy to verify that $F_1$ and $F_2$ are tensor functors. 
  \blmma
  \label{normal_functor}
  Suppose for some object $x$ of ${\rm Rep}(\hat{\clq})$,  $e_\clc \leq F_2(x)$. Then $x \in {\rm Rep}(\hat{\clq_H})$.
  \elmma
  \begin{proof} Decomposing  $x$ into irreducibles and using the fact that $F_2$ is a tensor functor, we can assume that there is some irreducible $\alpha \leq x$ such that $e_\clc \in F_2(\alpha)$. Let $\alpha$ correspond to $(\rho, \clh)$ where $\clh$ is a finite dimensional Hilbert space of dimension $d_\alpha$ and $\rho$ is unitarily conjugate to the projection $P_\alpha : \hat{\clq} \raro \hat{\clq}_\alpha \cong \clb(\IC^{d_\alpha}) $. The assumption $e_\clc \leq F_2(\alpha)$ means there is some one dimensional subspace, say spanned by a nonzero vector $\xi$, of $\clh$ such that \be \label{vvv} \rho(a)\xi=\epsilon(a)\xi\ee for $a \in \clc$. Recall the $\clc$-valued weight $T_H$ which 
   projects from $\hat\clq$ to $\clc$, given by $T_H(\cdot)=(\tau_H \ot {\rm id}) \circ \Delta(\cdot)$. Take some positive nonzero element $b \in \clq_h$ for some $h \in H$ and let $a=T_H(b) \in \clc$. Equation (\ref{vvv}) implies \be \label{uv} (\tau_H \ot \rho)(\Delta(b))\xi =\tau_H(b)\xi. \ee As $\tau_H(b)$ is nonzero, left hand side of 
   Equation (\ref{uv}) is also nonzero, which is possible only if there is some $h^\prime \in H$ such that $(P_{h^\prime} \ot P_\alpha )(\Delta(b)) $ is nonzero for all nonzero $b \in \clq_h$. But this means $h \leq h^\prime \ot \alpha$, hence $h \sim \alpha$ or $\alpha \sim h$, i.e. $\alpha \leq h \ot h_1$ for some $h_1 \in H$, so 
   in particular $\alpha \in H$. \end{proof}
     
     Let us define an equivalence relation $\sim^\prime$ on ${\rm Rep}(\clc)$ by the following
    \bdfn
     For $x, y \in {\rm Rep}(\clc)$, define $x \sim^\prime y$ if there are positive integers $m,n$ such that $x \leq my$ and $y \leq nx$, where for any object $z$ we denote the direct sum of $k$ copies of $z$ by $kz$. 
     
     \edfn
    \brmrk
     It is easy to verify that $\sim^\prime$ is an equivalence relation and $x \sim^\prime y$ if and only if $x$ and $y$ are direct sums of the same set of irreducibles, possibly with different multiplicities. In case $x,y$ are irreducible objects, $x \sim^\prime y$ if and only if they are isomorphic. 
     \ermrk

    \bcrlre
    \label{normal_functor_cor}
    $F_2(\alpha) \sim^\prime F_2(\beta)$ for some $\alpha, \beta \in {\rm Rep}(\hat{\clq}) $ if and only if $\alpha \sim \beta$.
     \ecrlre
  {\it Proof:}\\ We have $e_\clc \leq F_2(\alpha)\otimes \overline{F_2(\alpha)}  \leq F_2(\alpha)\otimes \overline{F_2(n \beta)}=nF_2(\alpha \ot \overline{\beta} )$ for some positive integer $n$. As $e_\clc$ is irreducible, this means $e_\clc \leq F_2(\alpha \ot \overline{\beta} )$, hence by Lemma \ref{normal_functor} we have $ \alpha \ot \overline{\beta} \in {\rm Rep} (\hat{\clq_H}) $, say $\alpha \ot \overline{ \beta} \leq h_1 \oplus  \ldots \oplus h_k$ for some $ h_1, \ldots, h_k \in H$.
   This implies, ${\rm Mor}( \alpha, \oplus_{i=1}^k h_i \ot \beta) \ne \{ 0 \}$, hence there is some $i$ such that ${\rm Mor}( \alpha  h_i \ot \beta) \ne \{ 0 \}$, i.e.$\alpha \sim \beta$. 
   
   For the converse, it is enough to note that $F_2(h)=e_\clc$ for all $h \in H$ . \qed\\
  
 % Note that the proof of the above corollary only uses the fact that the functor $F_2$ satisfies the condition in the statement of Lemma \ref{normal_functor}. Hence Corollary \ref{normal_functor_cor} is actually valid for any tensor functor $G$ satisfying the statement of Lemma \ref{normal_functor}. Let us record this observation for future use in the form of :
 % \bcrlre
 %$$ \label{normal_functor_cor_gen}
% Let $G$ be a tensor functor from ${\rm Rep}(\hat{\clq}) $ to some UTC $\clg$ such that $e \leq G(x)$ ($x \in {\rm Rep}(\hat{\clq}) $) implies $x \in {\rm Rep} (\hat{\clq_H}) $, where $e$ is the identity object in $\clg$. Define $\sim^\prime $ for objects in $\clg$ in the same way as done for ${\rm Rep}(\clc).$ 
 %%  Then $G(x) \sim^\prime G(y)$ for $x,y \in {\rm Rep}(\hat{\clq}) $ if and only if $x \sim y$.
 % \ecrlre

  For $\alpha \in I,$ let $$K_\alpha:=\{ x \in {\rm Irr}(\clc):~x \leq F_2(\alpha) \}.$$ 
  
  \blmma
  For any two $\alpha, \beta \in I$, the following are equivalent:\\
  (i) $K_\alpha \bigcap K_\beta$ is nonempty.\\
  (ii) $\alpha \sim \beta$.\\
  (iii) $K_\alpha=K_\beta$.
  \elmma
  {\it Proof:}\\
  (i) implies (ii): let $x \in K_\alpha \bigcap K_\beta$. Then $e_\clc \subseteq \overline{x} \ot x \leq \overline{F_2(\alpha)} \ot F_2(\beta)=F_2(\overline{\alpha} \ot \beta)$. By Lemma \ref{normal_functor}, $\overline{\alpha} \ot \beta \in {\rm Rep}(\hat{\clq}_H)$, hence $\alpha \sim \beta$.\\
  (ii) implies (iii) : follows from Corollary \ref{normal_functor_cor} by noting that $F_2(\alpha) \sim^\prime F_2(\beta)$ if and only if $K_\alpha=K_\beta$.\\
  (iii) implies (i) is obvious.\qed\\
%and $\q_{H}$ is a Woronowicz subalgebra of $(\q,\Delta)$.
\\
 As $F_2$ is surjective on objects, we get a partition of ${\rm Irr}(\clc)$ in terms of $K_\alpha$'s labelled by equivalence classes $[\alpha] \in I/ \sim$:
  $$ {\rm Irr}(\clc)=\bigcup_{ [\alpha] \in I/ \sim} K_\alpha.$$ 
 %Hence for any $x\in Q_{H}$,it's true that $\sum_{j}u^{\alpha}_{ij}xu^{*\alpha}_{i_{1}j}$ is in ${Q_{H}}$ and it's true for any arbitary irreducible object $u_{\alpha}$ of $C_{J}$. So,we can say that $(\tilde{Q_{0}})$ is a normal DQG subgroup of $(\tilde{Q_{J}})$ by using previous lemma.
 \end{proof}

 \section{Strongly projective corepresentation}

 Let $U$ be an object of the category $\alpha,$ considered in section \ref{3}.
 \bdfn
$U$ is said to be a strongly right projective corepresentation if $U\otimes V\otimes \bar{U}$ is a corepresentation of $\q$ for any corepresentation $V$ of $\q$. \edfn
\blmma If $U$ is a strongly right projective corepresentation then it is a right projective corepresentation of $\q.$\elmma
\begin{proof}
   By definition, $U\otimes V\otimes \bar{U}$ is a corepresentation of $\q$ for any corepresentation $V$ of $\q$. For $V=Id_{\bbc}\otimes 1_{\q}$,  $U\otimes \bar{U}$ is a corepresentation. So it is a right projective corepresentation of $\q$.
\end{proof}

%Let \begin{align*}
   % C_{stp}=\{\bigoplus x_{i}: s.t.  x_{i}\otimes y\otimes \bar{x_{i}}\in corep(Q) ,  \forall y\in corep(Q)\}
%\end{align*}
%\par
%we can  prove that $C_{stp} $ is a UTC where morphisom between two objects $u,v$ of this category defined by 
%\begin{align*}
  %  Mor(u,v)=\{T\in B(H_{u},H_{v}):v(T\otimes 1)=(T\otimes 1)u\}
%\end{align*}

\begin{rmrk}
 Let  $X,Y$ be two strongly right projective corepresentations of $\q$ then $X\otimes Y$ is  a right projective corepresentation of $\q$ because if we put  $V=Id_{c}\otimes 1_{\q}$ then $X\otimes Y\otimes  V \otimes \bar{Y}\otimes \bar{X}=X\otimes Y\otimes  \bar{Y}\otimes \bar{X}$ is a corepresentation of $\q$ .   
\end{rmrk} 

\blmma \label{comm} $U$ is a strongly right projective  corepresentation and $\Omega$ is a corresponding right 2-cocycle of $\q$ if and only if  right 2-cocycle $\Omega$ commutes with $\Delta(x)$ for any $x\in \q$.\elmma
\begin{proof}
    If $U$ is a right projective corepresentation and $(Id\otimes \Delta_{\Omega^{\ast}})(U)=U_{12}U_{13}$ then $\bar{U}$ is a left projective corepresentation and correspond to left 2-cocycle $\Omega.$
By definition it follows that for any corepresentation $V$ of $\q$, $U\otimes V\otimes \bar{U}$ is a corepresentation of $\q$.\\
So
\begin{align*}
    (Id\otimes Id\otimes Id\otimes \Delta)(U\otimes V\otimes \bar{U})=&(U\otimes V\otimes \bar{U})_{1234}(U\otimes V\otimes \bar{U})_{1235}\\
    =&(U_{14}V_{24}\bar{U_{34}}\otimes 1_{\q})(U_{15}V_{25}\bar{U_{35}})\\
    =&U_{14}U_{15}V_{24}V_{25}\bar{U_{34}}\bar{U_{35}}.
\end{align*}

and  
\begin{equation*}
     (Id\otimes Id\otimes Id\otimes \Delta)(U\otimes V\otimes \bar{U})\end{equation*}
        \begin{equation*}
           =(Id\otimes Id\otimes Id\otimes \Delta)(U_{14}V_{24}\bar{U_{34}}) 
        \end{equation*}  
     \begin{equation*}
          =(Id\otimes Id\otimes Id\otimes \Delta)(U_{14}V_{24}\bar{U_{34}})
     \end{equation*}
    \begin{equation*}
         =(Id\otimes Id\otimes Id\otimes \Delta)(U_{14})(Id\otimes Id\otimes Id\otimes \Delta)(V_{24})(Id\otimes Id\otimes Id\otimes \Delta)(\bar U_{34})
    \end{equation*}
   \begin{equation*}
       =U_{14}U_{15}(Id\otimes Id\otimes Id\otimes \Omega )V_{24}V_{25}(Id\otimes Id\otimes Id\otimes \Omega^{*})\bar{V_{34}}\bar{V_{35}}.
   \end{equation*}
  
 From those two equations, one can conclude that 
 \begin{align*}
     V_{24}V_{25}=&(Id\otimes Id\otimes Id\otimes \Omega )V_{24}V_{25}(Id\otimes Id\otimes Id\otimes \Omega^{*}).
    \end{align*}
Hence, 
\begin{align*}
    (Id\otimes \Omega)(Id\otimes \Delta)(V)=(Id\otimes \Delta)(V)(Id\otimes \Omega).
  \end{align*}
For any $x\in pol(\q),$ $\Omega \Delta(x)=\Delta(x)\Omega$.
  Converse part follows easily.
\end{proof}

  \blmma \label{5.5}
  If $\Omega$ is an invariant right 2-cocycle of $\q$ then\\
   1) $\Omega^{*}$ is a invariant right 2 cocycle of $\q$. \\
  2) If $U$ is a right projective corepresentation and $\Omega$ is corresponding 2-cocycle if and only if $U$ is a left projective corepresentation .\elmma
  \begin{proof}
      1) As $\Delta(x)$ commutes with $\Omega$ for all $x\in \q$,  $(\Delta \ot id)(y)$ commutes with $(\Omega\ot 1) $ for all $y\in \q\ot\q.$ Now taking $\ast$ on both side of the relation  $ (\Omega\otimes 1)(\Delta\otimes id)(\Omega)=(1\otimes \Omega)(1\otimes \Delta)(\Omega)$, is easily follows
       
    \begin{align*}
         (\Delta\otimes id)(\Omega^{\ast})(\Omega^{\ast} \otimes 1))=& ((\Delta\otimes id)(\Omega))^{\ast}(\Omega^{\ast} \otimes 1))\\
         =&((\Omega \ot 1)(\Delta \ot 1)(\Omega))^{\ast}\\
         =&((\Delta\otimes Id)(\Omega)(\Omega\ot 1))^{\ast}
        \\
\end{align*}
%and \begin{align*}
 %   (Id \ot \Delta)(\Omega^{\ast})(1\otimes \Omega^{\ast})=& (1\ot \Omega^{\ast})(Id\ot \Delta)(\Omega^{\ast})
%\end{align*}

  Similarly one can prove that $(1\otimes \Delta)(\Omega^{*})(1\otimes \Omega^{*})=(1\otimes \Omega^{*})(1\otimes \Delta)(\Omega^{*})$. So $\Omega^{*}$ is a right 2-Cocycle, which is invariant as $\Omega^{\ast}$ and $\Delta(x)$ commutes with $\forall x\in \q.$

 2) Note that by (1), $\Delta_{\Omega^{\ast}}(x)={_\Omega}^{*}\Delta(x).$ If $U$ is a right projective corepresentations then \begin{align*}
           (Id\otimes \Delta_{\Omega^{*}})(U)=(Id\otimes _{\Omega^{*}}\Delta)(U)=U_{12}U_{13}.
       \end{align*} Hence $U$ is also left projective corepresentation.
  \end{proof}
\bcrlre A unitary $\Omega\in \q\ot\q$ is an invariant right 2-cocycle if and only if it is an invariant left 2-cocycle.
\ecrlre
\begin{proof}
Clearly $\Omega$ is an invariant 2-cocycle, which commutes with $\Delta(x)$ for all $x\in \q$ if and only if $\Omega^{\ast}$ does so. Moreover by lemma \ref{5.5} we see that $\Omega$
invariant right 2-cocycle implies that $\Omega^{\ast}$ is an invariant right 2-cocycle,hence $\Omega$ is an invariant left 2-cocycle.

Similarly, we can prove that any left 2-cocycle $\Omega$ is also invariant right 2-cocycle.
\end{proof}
\brmrk
In view of the above, we will call a left/right invariant 2-cocycle simply an invariant 2-cocycle.
\ermrk

  \bdfn If $Y\in B(\clh_{y})\ot \q$ is a unitary element and for any corepresentation $V$ of $(\q,\Delta)$, $\bar{Y}\ot V\ot Y$ is a corepresentation of $\q$ then $Y$ is a strongly left projective corepresentation of $\q$.
  \edfn
\brmrk 
One can easily conclude that any strongly left projective corepresentation of $\q$ is a left projective corepresentation of $\q.$
\ermrk
\blmma
If $\Omega$ is a left 2-cocycle and $U$ is a strongly left projective corepresentation with corresponding left 2-cocycle $\Omega$ if and only if $\Omega$ commutes with $\Delta(x)$, for any $x\in q$.\elmma
\begin{proof}
 Proof of this lemma similar to the proof of the lemma (\ref{comm}) . Hence omitted. 
\end{proof}
  %Let $\Omega$ be a 2 cocycle of $(\q,\Delta)$ and it's an invariant 2 cocyle. Now, if we choose the right regular $\Omega^{*}$ corepresentation $V_{\q}\Omega^{*}$ then any right projective subrepresentation of $V_{\q}\Omega^{*}$ is a strongly left/right projective representation of $\q$. More generally, any right projective $\Omega$ corepresentation is a strongly right projective corepresentation.
  
 \blmma If $U$ is a unitary in $B(\clh_{u})\ot \q $ then the following are equivalent:\begin{enumerate}
     \item  $U$ is a strongly left projective corepresentation.
     \item $U$ is a strongly right projective corepresentation.
 \end{enumerate}
 \elmma
 \begin{proof}
     Let $U$ be a strongly right projective corepresentation with corresponding invariant 2-cocycle $\Omega_{u}.$  $\bar{U}$ is a left projective corepresentation with corresponding 2-cocycle $\Omega^{\ast}_{u}$. Hence by lemma \ref{5.5}, it is also a right projective corepresentation with corresponding 2-cocycle $\Omega^{*}_{u}.$ Hence, for every corepresentation $V$ of $\q,$ $\bar{U}\ot V\ot U$  corresponds to the cocycle $\Omega^{*}_{u}. \Omega_{u}=1$ (as $(Id\ot \Omega)$ commutes with $((Id\ot\Delta)(V)),$ that is $\bar{U}\ot V\ot U$ is a corepresentation of $\q.$ This proves $\bar{U}$ is a strongly right projective corepresentation, so $U$ is strongly left projective corepresentation of $\q.$

     The other way is similar to prove.
    %$\bar{U}\otimes V\otimes U$ is a corepresentation of $Q$. So, $\bar{U}$ is a strongly right projective corepresenation of $Q$ and it is a right,left projective $\Omega$ corepresentation of $Q$.
 %Similarly, we can prove the other direction.
 \end{proof}
In the view of above, we make the following definition

\bdfn A unitary $U\in B(\clh_{u})\otimes \q$  a strongly projective corepresentation if it satisfies any of two equivalent conditions:

\begin{enumerate}
    \item [1)] $U\otimes V\otimes\bar{U}$ is  corepresentation for all corepresentation $V$ of $\q$.
    \item[2)] $\bar{U}\otimes V\otimes {U}$ is  corepresentation for all corepresentation $V$ of $\q$.
\end{enumerate}
 
 \edfn

%More generally, we can say that a unitary $U\in B(H_{u})\ot\q$ is a strongly left/right projective representation if $U\ot V\ot\bar{U}$ is a corepresentation of $\q$ or we can say that If $U$ is a strongly projective corepresentation then these are the equivalent conditions:

 \blmma \label{mstp} Let  $U,W$ be two strongly projective $\Omega_{u},\Omega_{w}$ corepresentations,
 \begin{enumerate}
     \item [1)]  $U\otimes W$ is a strongly projective $\Omega_{u}\Omega_{w}$ corepresentations.
     \item[2)] $U$ is a subobject of $W$ then $\Omega_u=\Omega_w.$
 \end{enumerate}
 \elmma
 \begin{proof}
     1) It is already known that $(Id\otimes\Delta_{\Omega_{u}})(U)=U_{12}U_{13}$ and $(Id\otimes\Delta_{\Omega_{w}})(W)=W_{12}W_{13}$.
 Now,\begin{align*}
     (Id\otimes \Delta_{\Omega_{u}\Omega_{w}})(U\otimes W)=&(Id\otimes Id\otimes \Delta _{\Omega_{u}})(U_{13})(Id\otimes Id\otimes \Delta_{\Omega_{W}})(W_{23})\\
     &=U_{13}U_{14}W_{23}W_{24}\\
     &=U_{13}W_{23}U_{14}W_{24}\\
     &=(U\otimes W)_{12}(U\otimes W)_{13}.
 \end{align*} Hence, $(U\otimes W)$ is a strongly projective $\Omega_{u}\Omega_{w}$ corepresentation.

 2) Let $P$ is a projection in $Mor(W,W)$ such that $(P\ot 1)W=U=W(P\otimes 1).$
 \begin{align*}
     (Id\otimes \Delta_{\Omega^{\ast}_w})((P\ot 1)W)=& (P\ot 1)((Id\otimes \Delta_{\Omega^{\ast}_w})(W)\\
     =& (P\ot 1)W_{12}W_{13}\\
     =&((P\ot 1)W)_{12}((P\ot 1)W)_{13}\\
     =&U_{12}U_{13}\\
     =& (Id\otimes \Delta_{\Omega^{\ast}_u})(U).
 \end{align*} Hence $\Omega_{u}=\Omega_{w}.$
 \end{proof}
 
%Consider the subcategory 

%Let $F_{stp}$ is a Fiber functor from $C_{stp}$ to $Hilb_{f}$ defined by,
%\begin{align*}
   % &F(\bigoplus x_{i})=\bigoplus %H_{x_{i}}\\
   % &F:Mor(u,v)\rightarrow Mor(F(u),F(v))\\
   % &s.t F(T)=T
%\end{align*}, where $x_{i}$'s are strongly projective representation on $H_{i}$ and $T$ is in $Mor(u,v)$.\\

\brmrk When $\q=C(G)$ i.e. commutative as a $C^{\ast}$ algebra, any left/right projective corepresentation is automatically strongly projective corepresentation.For a compact group $G,$ there is a one-one correspondence between finite dimensional projective corepresentations of $(C(G),\Delta)$ and continuous projective representations of $G.$
Let $\phi$ be a continuous projective representation of $G$ over $V$ (Hilbert space) and $\omega_{\phi}$ be the corresponding 2-cocycle. We can write $\phi(g)=\sum c_{ij}(g)e_{ij},$ where $e_{ij}$ are the matrix units of $B(V)$. Then $u_{\phi}:=\sum e_{ij}\otimes c_{ij}$  is a unitary element of $B(V)\otimes \q$. One can easily conclude that $(Id\otimes \Delta_{\omega^{*}_{\phi}})(u_{\phi})=(u_{\phi})_{12}(u_{\phi})_{13}$. So, $\phi\mapsto u_{\phi}$ is a one-one correspondence between  projective corepresentations of $(C(G),\Delta)$ and continuous projective representations of $G$ on a finite dimensional Hilbert space .

\ermrk

 Recall the category $\alpha$ constructed in the Section \ref{3}, consider the subfamily of $obj(\alpha)$ given by $(U,\clh),$ where $U$ is a strongly projective corepresentation of $\q.$ The corresponding $C^{\ast}$ tensor category generated by this family , in a similar way as in Theorem \ref{main} , will be denoted by $C_{stp}$ and let $F_{stp}$ be the restriction of $F$ on $C_{stp}.$
Now, if we apply the Tannaka-Krein duality on $(C_{stp},F_{stp})$, we will get $\q_{stp}$ which is called the strongly projective envelope of $\q$. Clearly $\q\subseteq \q_{stp}$ is a Woronowicz subalgebra.
\par
Let $\clg$ be the UTC generated by corepresentation of $\q$ and $(Id_{\bbc}\otimes q)$ ,where $(Id_{\bbc}\otimes q)$ is an object of the UTC $C_{stp}$. 
\blmma If $(q,\bbc)$ is an object of this category $C_{stp}$, where $q$ is a unitary in $\q$ if and only if
     $(q^{*}\otimes q^{*})\Delta(q)$ commutes with $\Delta(x)$ for all $x\in \q.$ Moreover, in this case $\Delta(q)$ commutes with $(q^{*}\otimes q^{*})$.

\elmma

\begin{proof}
    $(Id_{\bbc}\otimes q)$ is an object of $C_{stp}$ that means that $(Id_{\bbc}\otimes q)V(Id_{\bbc}\otimes q^{*})$ is a corepresentation of $\q$ , for any corepresentation $V$ of $\q$. Hence    $(Id\otimes q)$ is a strongly projective corepresentation of $\q$ and $(q^{*}\otimes q^{*})\Delta(q)$ corresponding 2-coycle .
    We conclude from lemma \ref{comm} that $(q^{*}\ot q^{*})\Delta(q)$ commutes with $\Delta(x)$ for any $x\in \q.$
    %Now,one can easily conclude from previous lemma that $(q^{*}\otimes q^{*})\Delta(q)$ commutes with $\Delta(x)$ for any $x\in Q$. \\
  Taking $x\in \q$, we get
 \begin{align*}
     (q^{*}\otimes q^{*})=(q^{*}\otimes q^{*})\Delta(q)\Delta(q^{*})=\Delta(q^{*})(q^{*}\otimes q^{*})\Delta(q)
 \end{align*} Hence $\Delta(q)(q^{*}\otimes q^{*})=(q^{*}\otimes q^{*})\Delta(q)$.
\end{proof}

  Let $\q_{\clg}$ be the CQG obtained by Tannaka-Krein reconstruction by restricting $F$ to the subcategory $\clg$. We have the normalizer $N(C_{stp}, \clg)$ and let the CQG $N(\q_{\clg})$ be the one obtained by restriction of $F$ to $N(C_{stp},\clg).$ Here we recall the equivalence relation $\sim$ of Lemma \ref{4.11} on $Rep(\hat{\q})$ associated with the discrete normal subgroup $\hat{\q}_{\clg}$ of $\widehat{N(\q_{\clg})}$ and also the fiber functors \begin{align*}
      &F_{1}: Rep(\hat{\q}_{\clg})\to Rep(\widehat{N(\q_{\clg}}),\\
      &F_2: Rep(\widehat{N(\q_{\clg}}) \to Rep(\widehat{N(\q_{\clg})}/\hat{\q}_{\clg}).
  \end{align*} 

\blmma For an irreducible object $(U,\clh)\in obj(C_{stp}),$ the equivalence class $K_{U}$ contains exatly one (irreducibe) 1-dimensional representation, say $\rho_{U}$ of the quotient DQG $\widehat{N(\q_{\clg})}/\hat{\q}_{\clg}.$
%If $U_{\alpha}$ is an irreducible projective representation of $\q$, then  $K_{\alpha}=\{x\in irr({\hat{N({Q_{H}})}/\hat{Q_{H}})}:x\leq F(\pi_{\alpha})\}$ contains only a one dimensional representation of the DQG $\hat{N({Q_{H}})}/\hat{Q_{H}})$, where $\pi_{\alpha}(f)=(id\ot f)(U_{\alpha})$. 
\elmma
\begin{proof}

Denote by $\Pi_{U}: \widehat{N(\q_{\clg})} \to B(\clh)$ the representation dual to the corepresentation $U.$ Note that $U\otimes \bar{U}$ a corepresentation of $\q,$ so in particular an object of $\clg.$ Hence $F_2(\pi_{U}\otimes \bar{\pi}_{U})$ must be $n^{2}$ copies of the trivial representation of $e$ of $\widehat{N(\q_{\clg})}/\hat{\q}_{\clg},$ where $n= dim(\clh).$ Now suppose $F_{2}(\pi_{U})\cong \oplus^{k}_{i=1} \pi_i,$  each $\pi_i$ is a irreducible representation of $\widehat{N(\q_{\clg})}/\hat{\q}_{\clg}.$ As $F_{2}(\pi_{U})\otimes \overline{F_{2}(\pi_{U})}$ is direct sum of copies of $e$ only. But this is clearly possible if and only if each $\pi_{i}\otimes \bar{\pi_{j}}$ is isomorphic with $e$, hence in particular 1-dimensional. This further implies $\pi_i$ to be 1-dimensional and $\pi_i\cong\pi_j$ as well, which completes the proof.
\end{proof}
    %If $u_{\alpha}$ is an irreducible projective representation then $u_{\alpha}\otimes \overline{u_{\alpha}}$ is a corepresentation of $\q$. Now, We assume that $F(\pi_{\alpha})=(\pi_{\alpha})_{/({\hat{N({Q_{W}})}/\hat{Q_{W}}})}=\bigoplus{\pi_{i}}$, where $\pi_{\alpha}(f)=(id\otimes f)(x^{u_{\alpha}})$ for all $f\in \hat{N(Q_{W})}$ and $x^{u_{\alpha}}$ is the corresponding corepresentation of $u_{\alpha}$ in $N(Q_{W})$ and $\pi_{i}$ are irreducible subrepresentations of $F(\pi_{\alpha})$. we know $F(\pi_{\alpha}\otimes \overline{\pi_{\alpha}})=e_{\tilde{N({Q_{W}})}/\tilde{Q_{W}})}$ which implies $\bigoplus(\pi_{i}\otimes \overline \pi_{j})=e_{\hat{N({Q_{W}})}/\hat{Q_{W}})}$. Hence, we can say that for any $i,j$ $\pi_{i}\otimes \overline{\pi_{j}}$ is isomorphic to $e_{\hat{N({Q_{W}})}/\hat{Q_{W}})}$. Wlog, we assume that $\pi_{i}\otimes \overline{\pi_{j}}=e_{\hat{N({Q_{W}})}/\hat{Q_{W}})}$. So, $\pi_{i}$ is dual to the  unitary one-dimensional corepresenation of the dual CQG of ${\hat{N({Q_{W}})}/\hat{Q_{W}})}$ or more clearly we can say that there exists a group like element $q_{i}$ in the dual CQG of ${\hat{N({Q_{W}})}/\hat{Q_{W}})}$ for which $\pi_{i}(f)=(id\otimes f)(q_{i})$, for all $f$ in ${\hat{N({Q_{W}})}/\hat{Q_{W}})}$. So, we can conclude that $\pi_{i}=\pi_{j}$ as both are one-dimensional representations and $\pi_{i}\otimes \overline{\pi_{j}}=e_{\hat{N({Q_{W}})}/\hat{Q_{W}})}$. Hence $K_{\alpha}$ contains only a one-dimensional representation.

\blmma \label{5.17}
If two irreducible strongly projective corepresentations $(U,\clh),(V,\clk)$, we have $K_{U}=K_{W}$ if and only if the corresponding invariant 2-cocycles $\Omega_{U},\Omega_{W}$ are related by :
\begin{align}
\Omega_{w}=\Delta(q)\Omega_{u}(q^{\ast}\otimes q^{\ast})
\end{align} for some unitary $q$ such that $(q,\bbc)\in obj(C_{stp}).$
\elmma
\begin{proof}
If $K_{U}=K_{W}$ if and only if $U\sim W,$
which is equivalent to $U$ being a subobject of $W\otimes V$, where $V$ is an object of $C_{stp}.$ By lemma \ref{mstp} cocycle of $W\otimes V$ is $\Omega_{w}\Omega_{v}.$ We can replace $V\in obj(C_{stp})$ by any suitable generating subset . In particular $V$ may be either a corepresentation or a unitary element $q$ such that $(q,\bbc)\in obj(C_{stp}).$ If $V$ is a corepresentation then $\Omega_{u}=\Omega_{w}$ and the other case $\Omega_{u}=\Omega{w}(q\otimes q)\Delta(q^{\ast}).$

Conversely, suppose $\Omega_{w}=\Delta(q)\Omega_{u}(q^{\ast}\otimes q^{\ast}).$ For $U$ and $W(Id_{\bbc}\otimes q)$ they have the same invariant 2-cocycle $\Omega_{w}.$
 Hence $V=\bar{W}\otimes U(Id_{\bbc}\otimes q)$ is a corepresentation of $\q.$ Now  $U(id_{\bbc}\otimes q)\leq W\otimes \bar{W}\otimes U(id_{\bbc}\otimes q)= W\otimes V$ where $\leq$ means subobject. Hence $U$ is a subobject of $W\otimes V\otimes (Id_{\bbc}\otimes q)$ which completes the proof.
 \end{proof}

\bthm\label{5.18}
     $\widehat{N({Q_{\clg}})}/\hat{Q}_{\clg}$ is isomorphic to $C_{0}(\Gamma_{\q})$ for some discrete group $\Gamma_{\q}.$ 
\ethm
\begin{proof}

It follows from lemma \ref{5.17} that every irreducible representation of the DQG $\widehat{N({Q_{\clg}})}/\hat{Q}_{\clg},$ is 1-dimensional . Hence it is a commutative as a $C^{\ast}$ algebra ,which proves the theorem.
\end{proof}
  %  We already know from lemma (\ref{mstp}) that tensoring of two strongly projective representation is a strongly projective corepresentation. Let $U,V$ are two irreducible strongly projective corepresentations. 
 %   Let,\begin{align*}
     %   &\pi_{u}(f)=(id\ot f)(U)\\
    %     &\pi_{v}(f)=(id\ot f)(V), f\in \hat{N(Q_{H})} \end{align*}
         
%It is already known from previous lemma that $F(\pi_{u}\ot \pi_{v})\cong F(\pi_{u})\ot F(\pi_{v})$ is dual to the one dimensional unitary corepresentation of the dual CQG of $\Hat{N({Q_{H}})}/\Hat{Q_{H}}$. As dual CQG$\Hat{N({Q_{H}})}/\Hat{Q_{H}}$  generated by the one dimensional unitary corepresentation. So, the universal dual CQG $\Hat{N({Q_{H}})}/\Hat{Q_{H}}$ is a cocommutative CQG. Hence it is $C^{*}(G)$, for some group $G$. Without loss of generality we can assume, that $G=\{[K_{\alpha}]\}$, $\alpha\in I$, where $I$ is the index of the collection of all inequivalent irreducible objects. So, the DQG $\Hat{N({Q_{H}})}/\Hat{Q_{H}}$ is isomorphic to the dual of CQG $C^{*}(G)$.

  % \end{proof}
  \bthm For   a compact group $G$ with $\q=C(G)$, the group $\Gamma_{\q}$ of theorem is isomorphic with the group $H^{2}_{uinv}(C(G),S^{1}).$
  \ethm

  \begin{proof}
Let $(U,\clh),(Id_{\bbc}\otimes q)$ be two strongly projective corepresentations of $C(G).$  Then $U\otimes (Id_{\bbc}\otimes q)\otimes \bar{U}$ is isomorphic with $U\otimes \bar{U}\otimes (Id_{\bbc}\otimes q)\in obj(C_{\clg})$. From this we can conclude that $C_{stp}$ is a subcategory of $N(C_{stp},C_{\clg})$ and by our construction of normalizer $N(C_{stp},C_{\clg})$ is a subcategory of $C_{stp}$. Hence $C_{stp}=N(C_{stp},C_{\clg})$ that implies $\widehat{N(\q_{\clg})}=\hat{\q}_{stp}.$ We already observed that the set of 1-dimensional irreducible representation of $C_{0}(\Gamma_{\q}),$ which are the points of $\Gamma_{\q}$, are in 1-1 correspondence with the equivalence classes $K_{U}.$ Moreover by lemma \ref{5.17}, $K_{U}=K_{W}$ if and only if $\Omega_{u} ~and ~\Omega_{w}$ represents the same classes at $H^{2}_{uinv}(C(G),s^{1}).$ As $C(G)$  is a commutative, it is equal to its center, so that every $q\in C(G)$ is central and $(q\ot q)\Delta(q^{\ast})$ is an invariant 2-cocycle. Thus, there is a 1-1 correspondence between $H^{2}_{uinv}(C(G),S^{1}) ~and ~\rho_{U}$ which is clearly also a group homomorphism, as $\Omega_{u\otimes w}=\Omega_{u} \Omega_{w}.$
\end{proof}
\brmrk For a general CQG, which is not commutative as a $C^{\ast}$-algebra it is interesting to ask how the group $\Gamma_{\q}$ obtained in theorem \ref{5.18} is related to $H^{2}_{uinv}(\q,S^1).$ It is also interesting to explore the functorial property of the association $\q \rightarrow \Gamma_{\q}$ from the category of CQGs to the category of discrete groups.

\ermrk

\section{Examples}
In this section, we will calculate the projective envelope and the discrete group $\Gamma_{\q}$ of the theorem \ref{5.18} for certain compact quantum groups, including  the group ring of $Z_8\rtimes Aut(Z_8).$ Our computation of $\Gamma_{\q}$ for the group ring $Z_8\rtimes Aut(Z_8)$ demonstrates that $\Gamma_{\q}$ is different from the group $H^2_{uinv}(C^{\ast}(Z_8\rtimes Aut(Z_8)), S^{1}).$
 Note that for compact group $G,$ left/right/bi-projective corepresentations of $C(G)$ coincides as the $C^{\ast}$ algebra is commutative. Hence $\widetilde{C(G)}_{L}=\widetilde{C(G)}_{R}=\widetilde{C(G)}_{bi}$ and we denote it simply by $\widetilde{C(G)},$ called the projective envelope.
\blmma\label{comp} For a classical compact group $G$, the projective envelope $\widetilde{C(G)}$ of $C(G)$ is a commutative $C^{\ast}$ algebra, that is isomorphic with $C(\tilde{G})$ for some compact group $\tilde{G}$.\elmma
\begin{proof}
Let $U\in B(H_{u})\otimes C(G),V\in B(H_{v})\otimes C(G)$ be two irreducible inequivalent projective corepresentations of 
$C(G)$. Without loss of generality ,we assume that $U=\sum_{i,j} e_{ij}\otimes u_{ij},V=\sum_{k,l}f_{kl}\otimes v_{kl}$. Let $T:H_{u}\otimes H_{v}\rightarrow H_{v}\otimes H_{u}$ is the flip operator. For $e_{i_{0}}\in H_{u},f_{k_{0}}\in H_{v}$ 
\begin{align*}
   (T\otimes 1)(U\otimes V)(e_{i_{0}}\otimes f_{k_{0}})=&(T\otimes 1)(\sum_{i,j,k,l}e_{ij}(e_{i_{0}})\otimes f_{kl}(f_{k_{0}})\otimes u_{ij}v_{kl})\\=&(T\otimes 1)(\sum_{i,k} e_{i}\otimes f_{k}\otimes u_{ii_{0}}v_{kk_{0}})\\=&\sum_{i,k} f_{k}\otimes e_{i}\otimes u_{ii_{0}}v_{kk_{0}}.
\end{align*} and 
\begin{align*}
    (V\otimes U)(T\otimes 1)(e_{i_{0}}\otimes f_{k_{0}}\otimes 1)=&\sum_{i,j,k,l}f_{kl}(f_{k_{0}})\otimes e_{ij}(e_{i_{0}})\otimes v_{kl}u_{ij}\\=&\sum_{k,i} f_{k}\otimes e_{i}\otimes v_{kk_{0}}u_{ii_{0}}\\
    =&\sum_{i,k} f_{k}\otimes e_{i}\otimes u_{ii_{0}}v_{kk_{0}}.
\end{align*}\\
So, the unitary operator $T\in Mor(U\otimes V,V\otimes U)$. From Tannaka-Krein duality, we can say that $T\in Mor(X^{U\otimes V},X^{V\otimes U})$. As $F_{2}(U\otimes V)$ is the identity operator that's why $X^{U\otimes V}=X^{U}\otimes X^{V}$. Now, we can easily conclude that $x^{u}_{ii_{0}}x^{v}_{kk_{0}}=x^{v}_{kk_{0}}x^{u}_{ii_{0}}$, where $i_{0},k_{0}$,$u,v $ are arbitrary. As matrix coefficients of inequivalent irreducible corepresentations is a basis of the Hopf $*$-algebra $\widetilde{C(G)}$, we conclude that $\widetilde{C(G)}$ is a commutative compact quantum group.
\end{proof}

Let us recall a few well-known facts about the corepresentation theory of $SU_q(2)$ and $SO_q(3)$ for $q \in (0,1].$. Recall from \cite{ woronowicz1987compact}, \cite{woro}  and references therein that the set of finite dimensional, irreducible unitary corepresentatons of $SU_q(2) $ is labelled by $\{ \frac{n}{2},~n=0,1, 2, \ldots \}$, where for any half-integer $\frac{n}{2}$, 
there is exactly one (upto isomorphism) irreducible  corepresentation, say  $U_n$, of dimension $2n+1$. $SO_q(3)$ is a Woronowicz subalgebra of $SU_q(2)$ and ${\rm Corep}(SO_q(3) ) $ is generated by $U_n, n=0, 1, 2, \ldots $. Moreover, it is known that $SU_q(2) $ is projective torsion free (\cite{proj2DeMaNe}, \cite{voigt} ), hence any finite dimensional unitary right projective corepresentation of $SU_q(2)$ is isomorphic with $U (1 \ot v)$ for some unitary corepresentation $U$ and unitary element $v$ of $SU_q(2)$.
\bthm
\label{soq3} 
Every finite dimensional   unitary right projective corepresentation of $SO_q(3)$ is isomorphic with $U (1 \ot u)$ for some unitary corepresentation $U$ of $SU_q(2)$ and unitary element $u$ of $SO_q(3)$.
\ethm
\begin{proof}
As $SO_q(3)$ is a subalgebra of $SU_q(2)$, any  (unitary, finite dimensional) right  projective corepresentation $V$ is also a unitary, finite dimensional right projective corepresentation of $SU_q(2)$, so by the discussion preceeding the theorem, it must be isomorphic with $U(1 \ot v)$ for some unitary corepresentation $V$ of $SU_q(2)$ and unitary $v \in SU_q(2)$. However, this implies, ${\rm ad}_V$ is isomorphic with ${\rm ad}_U$. Now, it also follows from the detailed fusion rules of ${\rm Corep}(SU_q(2)) $ (\cite{woro}) that for any $n \in \{ 0 \} \bigcup \frac{\IN}{2},$ 
 $U_n^{(c)} \cong U_n$ and $U_n \otimes U_n$ is a direct sum of $U_m$ for $m \in \IN \bigcup \{ 0 \}$l i.e. belongs to ${\rm Corep} (SO_q(3) )$.  Writing $U$ as a direct sum of irreducibles, it follows that $U \otimes U^{(c)} $ is in ${\rm Corep} (SO_q(3)),$ so that the matrix coefficients of ${\rm ad}_U$ are actually in $SO_q(3)$. In other words, we have that the coactions ${\rm ad}_V$ and ${\rm ad}_U$ of $SO_q(3)$ are isomorphic, implying that $V \cong U(1 \ot u)$ for some $u \in SO_q(3)$. \end{proof}
 
Let us denote the (nontrivial)  right $2$-cocycle corresponding to $U_{\frac{1}{2}} $ by $\Omega_{\frac{1}{2}}$. 
 
 \bcrlre
 \label{gamma}
 For any irreducible unitary right projective corepresentation $V$ of $SO_q(3)$, the corresponding right $2$-cocycle is either isomorohic to the trivial cocycle or the cocycle $\Omega_{\frac{1}{2}}$. In other words, 
  there are exactly $2$ right $2$-cycles, $1$ and $\Omega_{\frac{1}{2}}$ upto isomorphism. A similar statement holds for left projective corepresentations as well. 
 \ecrlre
\begin{proof}
It follows from \ref{soq3} that $V  \cong   U_n (1 \ot u)$ for some half-integer $n$ and unitary $u$, hence the $2$-cycle corresponding to $V$ is isomorphic with that of $U_n$. Now,  $U_n$ is a corepresentation, i.e. has trivial cocycle, for whole integer $n$. 
  For fractional $n \geq \frac{3}{2}$,  the fusion rules imply that $U_n$ is a sub-object of $U_{\frac{1}{2}}  \otimes U_{n-\frac{1}{2} } $ and $U_{n-\frac{1}{2}} $ is a corepresentation. This implies that the cocycle corresponding to $U_n$ is isomorphic with that of $U_{\frac{1}{2}}$, i.e. $\Omega_{\frac{1}{2}}$. 

The case of left projective corepresentations can be proved in a similar way, hence omitted. 
 \end{proof} 
\blmma
\label{trivial}
Let $0<q<1$. If a unitary element $w \in SO_q(3)$ is a one dimensional strongly projective corepresentation, then $w$ must be a scalar multiple of $1$.
\elmma
\begin{proof}
By definition of strongly projective corepresentation, $(1 \ot w)U(1 \ot w^*) $ is a corepresentation for every unitary corepresentation $U$ of $SO_q(3)$, and clearly, it is irreducible if and only if $U$ is so. As there is only one unitary irreducible upto isomorphism for each integer $n$, we have 
 $(1 \ot w)U_n(1 \ot w^*) \cong U_n$. i.e. $(1 \ot w)U_n(1 \ot w^*)=(\gamma \ot w)U_n(\gamma^* \ot w^*)  $ for some unitary $\gamma$ in $\clb(\IC^{2n+1} ) $. Taking trace in the first leg, we get $w \chi_n w^*=\chi_n$, i.e. $w$ commutes with $\chi_n$ denotes the character of $U_n$, i.e. $\sum_i U_n(i,i), $ when 
  we write $U_n$ as a matrix with $SO_q(3)$ valued entries $(( U_n(i,j) ))$ w.r.t. any orthonormal basis of $\IC^{2n+1}$. Moreover, for any two unitary corepresentations $U,V$ of $SO_q(3)$, it is clear that ${\rm Mor} (U_w, V_w)={\rm Mor}(U,V)$, where we have denoted $(1 \ot w)U(1 \ot w^*) $ by 
  $U_w$ for any corepresentation $U$. Using this observation for the special case of $U$ to be the one-dimensional trivial corepresentation and $V$ to be $U_n \ot U_n $ for any half-integer $n$ ( note that, it follows from the fusion rules that $U_n \ot U_n$ is direct sum of $U_m$ with whole integer $m$ only), we  conclude that the spectral projection correponding to the trivial coreprsentation in $V=U_n \ot U_n$ and $V_w$ coincide. Using the expression of this projection (see, for example, \cite{pod} ), this is equivalent to $ ({\rm id } \ot h) (V_w)=({\rm id } \ot h) (V)$, where $h$ denotes the Haar state. Decomposing $V$ into irreuducibles and noting that varying $n$, this will give us all the irreducibles of $SO_q(3)$, we conclude that $h( wU_k(i,j) w^*)=h(U_k(i,j)) $ for all whole integer $k$ and $i,j$. As $U_k(i,j) $'s span a dense subset of $SO_q(3)$, it follows that $h(w a w^*)=h(a)$ for all $a \in SO_q(3)$. Now, viewing $SO_q(3)$ as a subalgebra of $SU_q(2)$ and recalling the modular operator $\rho$ for the Haar state, i.e. $h(ab)=h(\rho(b)a)$ for all $a,b$, where $\rho(\cdot)=f_{-1} \ast \cdot f_{-1} $ in the notation of  \cite{woronowicz1987compact},  we get, $h(\rho(w^*) w a)=h(a)$ for all $a$, hence by faithfulness of $h$, $\rho(w^*)w=1=w^*w$, i.e. $\rho(w^*)=w^*$, hence $\rho(w)=w$ too. Using the notation of \cite{ woronowicz1987compact}, we can write any element of $SU_q(2)$
   as a strongly convergent sum (in the GNS space of the Haar state) of the form $\sum_k \alpha^k \phi_k(\gamma)+ \sum_l (\alpha^*)^l \psi_l(\gamma)$, where $\alpha, \gamma$ are the standard generators of $SU_q(2)$ as in \cite{pod} and $\phi_k, \psi_l$ are continuous functions from the spectrum of the normal element $\gamma$ to $\IC$. Using the expression of $\rho$ on the monomials in $\alpha, \alpha^*, \gamma$, it is easy to see that the fixed point subalgebra of $\rho$ consists of functions of $\gamma$. Hence we can write $w$ as a function of $\gamma$, say $f(\gamma)$, for some continuous $f : \sigma(\gamma) \raro \IC$. 
   From the commutation relations of the generators we get $\alpha f(\gamma)=f(q \gamma) \alpha$. But then, $w=f(\gamma)$ commutes with  the character of $U_1$, which is given by (see, for example, \cite{woronowicz1987compact}, $\alpha^2+(\alpha^*)^2)$ plus some function of $\gamma$, hence $f(q^2\gamma)=1$ But in a faithful representation (see \cite{ woronowicz1987compact}) 
    of $SU_q(2)$ into $\clb(l^2(\IN \ot \IZ)$, $\gamma$ is mapped to the operator $q^N \ot S$, where $N$ denotes the `number operator' and $S$ the bilateral shift. Thus, the spectrum of $\gamma$ is $\{ 0 \} \bigcup \{ q^nz:~n \geq 1, z \in S^1 \}$ and the condition $f(q^2\gamma)=f(\gamma)$ implies,      $f(q^kz)=f(q^{k+2}z)$ for all $n$, hence $f(q^nz)=f(q^{n+2m} z)$ for all $m\ geq 1$, and letting $m \raro \infty$, $f(q^nz)=f(0)$, i.e $f$ is a constant function.
 \end{proof}

\bthm 
\label{strprojchoice}
$U_n$ is strongly projective corepresentation for all half integer $n$. Moreover, $\Gamma_{SO_q(3)} \cong \IZ_2$, with the nontrivial element represented by the equivalence class of  $U_{\frac{1}{2}}.$ 
\ethm
\begin{proof}
From the fusion rules of $SU_q(2)$, it is easy to see that $U_n $ is self-conjugate and for any half- integer $n$ and $U_n   \otimes U_m \otimes  U_n $ is the direct sum of irreducibles of $SO_q(3)$ only, for any whole integer $m$. This proves strong projectivity. 

By Lemma \ref{trivial}, any strongly projective one dimensional unitary corepresentation $w$  is of the form $c1$ for constant $c$, hence $U_n \ot w \ot U_n \cong U_n \ot U_n$, which is in the category of corepresentations of $SO_q(3)$. This, combined with the first part of the present proof, we see that each $U_n$ 
 is in the normalizer used for defining $\Gamma_{SO_q(3)} $.  Moreover, for any half-integer $n > \frac{1}{2}$, $U_n$ is isomorphic with a direct summand in $U_{\frac{1}{2}} \ot U_m$ for a suitable whole integer $m$, which implies that the equivalence class of $U_n$ w.r.t. the equivalence relation $\sim$ discussed before in the context of defining $\Gamma_\clq$, is same as that of $U_{\frac{1}{2}}$. This completes the proof.
 \end{proof}

 It follows from the above proof that  for any finite dimensional  right (hence also left) projective unitary coreprsentation $V$ of $SO_q(3)$, there is some strongly projective corepresentation $W$ such that $V \sim W$. Indeed, $W$ can be chosen as direct sum of copies of $1$ and $V_{\frac{1}{2}}$. 
 
 Let us make this property of $SO_q(3)$ into a formal definition. Recall the UTC $\tilde{\clc}:=\tilde{\clc^L_\alpha}= \tilde{\clc^R_\alpha}$,  $\clc_{\rm stp}$ corresponding to left/ right and strongly projective corepresentations respectively.   Note that $\clc_{\rm stp}$ is a full subcategory of $\tilde{\clc}$.  Let us also consider another full sub-category 
  of $\tilde{\clc}$  denoted by $\clh$, which is generated by the corepresentations of $\clq$ and unitary elements $q$ of $\clq$, viewed as one-dimensional bi-projective corepresentations. 
 \bdfn
 We call a CQG $\q$ to have enough strongly projective corepresentations if  every object $v $ of $\tilde{\clc}$ is a sub-object of $h \ot u$ for some objects $h$ of $\clh$ and $u$ of $\clc_{\rm stp}$
 \edfn

It is easy to see that if $v$ in the above definition is irreducible, we can choose both $h$ and $u$ to be irreducible too and then the class of $u$ in the group $\Gamma_\clq$, say $[u]$,  is uniquely determined. In other words, $v \mapsto [u] $ gives a well-defined map from irreducible objects of $\tilde{\clc}$ (hence extends to 
 the whole of the ${\rm Obj}(\tilde{\clc}) $ to $\Gamma_\clq$. Let us denote this map by $\pi_\Gamma$. In case of $SO_q(3)$, it follows from our previous discussion that $\pi_\gamma$ maps $U_{\frac{n}{2}} $ to $0$ for $n$ even and to  $1$ for $n$ odd. 
Moreover, note that the condition $v \leq h \ot u$ is interchangeable with $u \leq v \ot h^\prime$ ($\leq$ means sub-object) for some object $h^\prime$ of $\clh$ and also by taking conjugates, these conditions are equivalent to the conditions of the form $u \leq h_1 \ot v$ or $v \leq h_2 \ot u$ for some $h_1, h_2 $ in ${\rm Obj}(\clh)$. 

\bthm 
$\widetilde{C(SO(3))}\cong C(SU(2)\times \widehat{U}(C(SO(3)))),$ where $\widehat{U}(C(SO(3))$ is  dual of the discrete unitary group of $C(SO(3)).$ 
\ethm
\begin{proof}
It is already known from Lemma \ref{comp} that $\widetilde{C(SO(3))}$ is a commutative compact qunatum group.
%For an irreducible projective corepresentation $U$ of $C(SO(3))$, if $dim(U)\geq 2$ and corresponding 2-cocycle $\Omega$ is not coboundary equivalent to trivial 2-cocycle then $U$ is an irreducible corepresentation of $C(SU(2)).$ Let $U$ be an irreducible projective corepresentation of $C($. 
Restricting Theorem \ref{soq3} for $q=1$, we see that  any unitary projective corepresentation $V$  can be written as a $V=UId\ot u),$ where $u$ is a unitary element of $C(SO(3)).$ and $U$ is a unitary corepresentation. of $SU(2)$. It follows from Theorem \ref{main} that the projective envelope of $C(SO(3))$ is generated by the matrix coefficients of $U_{n}$'s (over all half integer $n$) and the one dimensional corepresentations $Id_{\bbc}\ot X^u$, where $u\in U(C(SO(3))).$  Matrix coefficients of $U_{n}$ generated the Hopf *-algebra $Pol(SU(2)).$  $\{Id_{\bbc}\ot X^u: u\in U(C(SO(3)))\}$ generated the Hopf *-algebra $Pol(C^* (U(SO(3))),$ where $U(C(SO(3)))$ considered as a discrete group. If we take the universal completion of this Hopf * -algebra which is generated by the matrix coefficients of $U_n$ and $Id_{\bbc}\ot X^u$ then we will get the compact quantum group $C(SU(2)\ot_{max} C^{*}(U(C(SO(3))))\cong C(SU(2))\ot_{min}C(\widehat{U}(C(SO(3))))\cong C(SU(2)\times \widehat{U}(C(SO(3)))).$
\end{proof}

Next, we compute the discrete group $\Gamma_{\q}$ of the Theorem \ref{5.18} for the group ring of $Z_{8}\rtimes Aut(Z_8).$
The group $G:=Z_8\rtimes Aut(Z_8)$ was  considered by G. E. Wall in his paper \cite{holg}, $Z_8\rtimes Aut(Z_8)$ is generated by $s,t,u$, where $s,t,u$ satisfies the relations 
\begin{align}
s^2=t^2=u^8=1, ~~~st=ts,~~sus^{-1}=u^3,~~~tut^{-1}=u^5~.
\end{align}
It has been proven in \cite{cohomom} that $H^2_{uinv}(C^{\ast}(Z_{8}\rtimes Aut(Z_8), S^1)=Z_2.$ It is also proven in \cite {dtwist} that the
nontrivial invariant 2-cocycle $\Omega $ is given by $\Omega=(v\ot v)\Delta(v^{\ast}),$  where
\begin{align*}
 v=\frac{1}{2}(1+u^4)+\frac{\sqrt2}{4}u(1-u^2-u^4+u^6).
\end{align*}
One can check that $v$ is a self adjoint unitary element and $v^2=1.$
\bthm For the group ring $Z_8\rtimes Aut(Z_8),$ $\Gamma_{\q}$ is trivial.
\ethm
\begin{proof}
It follows from lemma \ref{5.17} that the equivalence class $K_{Id_{\bbc}\ot 1}=K_{Id_{\bbc}\ot v}.$ By our construction, the group $\Gamma_{\q}$ is in one-one correspondence with the equivalence classes $K_{U},$ where $U$  irreducible strongly projective corepresentations. Hence $\Gamma_{\q}$ is the trivial group.
\end{proof}
\section{Application to physics: 1 D topological phases}
In \cite{Phy1} 
classification of gapped symmetric topological phases in 1-D spin systems has been done using projective representation of the symmetry group say ($G$). The second cohomology $H^2(G, \bbc)$ has been used to label the topological phases.\\
 We propose an almost verbatim generalization of the framework, replacing the symmetry group $G$ by suitable compact quantum groups, say $\q.$ Let us consider matrix product states (MPS) as in \cite{Phy1}
 with on-site symmetry given by some unitary corepresentation $U^{0}$ of $\q.$ Thus, the product states have the quantum symmetry given by $U^{0}\ot..\ot U^{0}.$ We can easily adapt the arguments of Appendix D of 
  \cite{Phy1}
 to show that, in the notation of  \cite{Phy1},
 a state which does not break the on-site symmetry is associated with a sequence of unitaries $M_{[k]},N_{[k]}, k\geq 1,$ on finite dimensional Hilbert space $\clh_{[k]}$, say, such that for all $k,$ $U_{[k]}:=N^{-1}_{[k]}\ot M_{[k]},$ gives a corepresentation upto the $\q-$ valued phase, meaning $U_{[k]}=W_{[k]}(1\ot q),$ for some unitary corepresentation $W_{[k]}$ on $\clh_{[k]}\ot \clh_{[k]},$ or in other words $Ad_{U_{[k]}}$ is a coaction on $B(\clh_{[k]}\ot \clh_{[k]}).$
 \blmma
 \label{888}
  Suppose that $V\in B(\clh)\ot \q,T\in B(\clk)\ot \q $ unitaries ($\clh,\clk$ are finite dimensional Hilbert spaces) such that $V \ot T=W(1 \ot q)$ for some unitary corepresentation $W$ and some unitary element $q$ of $\clq$.\\
  (i) Then 
  $V$ and $T$ are  right projective corepresentations of $\q$ with associatied 2-cocycles $\Omega_{V},\Omega_{T} $ respectively (say).\\
  (ii) If $V,T$ are strongly projective, we must have $\Omega_V=\Omega_T^{-1}$.\\
  (ii) Assume furthermore that  $\clq$ has enough strongly projective corepresentations. Then the images $\pi_\Gamma(v)$ and $\pi_\Gamma(T)$ of $V$ and $T$ under the map $\pi_\Gamma$ defined in Section 7, will be inverses to each other, i.e. $\pi_\Gamma(V)=\pi_\Gamma(T)^{-1}$. 
 \elmma
\begin{proof}
    (i) Consider $Ad_{V\ot T}(a\ot 1_{\clk}\ot 1_{\q})=(Ad_{V})_{13}$, where $a\in B(\clh).$ It is an coaction. Hence $V$ must be right projective corepresentation. Similarly, we can conclude that $T$ is right projective corepresentation as well. 
    
    (ii) This follows as $V \ot T$ is a strongly projective corepresentation with the cocycle $\Omega_V \Omega_T$ and $W(1 \ot q)$ corresponds to the trivial cocycle.
    
    (iii)  First assume $V,T$ to be irreducible and choose irreuducible, strongly projective $V_1, T_1$  and objects $h_1,h_2$ of $\tilde{\clc}$ such that $V_1 \leq h_1 \ot V,$ $ T_1 \leq T \ot h_2$. But then, $V_1 \ot V_2 \leq h_1 \ot V \ot T \ot h_2,$ which is by assumption an object in $\tilde{\clc}$. This clearly implies,
     $V_1^{(c)} \sim T_1$, which proves (iii). 
    \end{proof}
 %We apply the above to the set up discussed before to get that $N^{-1}_{[k]}$ and $M_{[k]}$ are bi-projective cotrepresentation.  
%In case  $N^{-1}_{[k]}$ and $M_{[k]},$  are strongly projective,  the corresponding cocycles  $(\Omega^{N}_{[k]})^{-1}$ and $\Omega^{M}_{[k]},$ (say) are related by $\Omega^{M}_{[k]}=\Delta(q)\Omega^{N}_{[k]}(q^{-1}\ot q^{-1})$ for some $q.$
This can be summarized as:
\bthm Under the assumption that $\clq$ has enough  strongly projective corepresentations,  the elements of the group $\Gamma_\clq$ corresponding to $M_{[k]}$ and $N_{[k]}$ are equal. 
\ethm
This shows that elements of $\Gamma_\clq$ labels different topological phases. Clearly, we can apply this to  the case of $SO_q(3)$, which has enough strongly projective corepresentations as shown earlier. In this case, there are two possible phases as $\Gamma_\clq=Z_2$. 
\brmrk
Note that the above definition of two invariant cocycles being cohomologus is exactly the same which we have used for our construction of $\Gamma_{\q},$ unlike the definition of $H^2_{uinv}(\q, S^1)$ in \cite{snbook}. Indeed, two invariant cocycles $\Omega_{1}, \Omega_{2}$ related by $\Omega_1=\Delta(q)\Omega_2 (q^{-1}\ot q^{-1})$ for any $q$ (not necessarily from the center of $\q$). Clearly gives the same symmetry-protected phase in this physics model. In other words, the classes of strongly projective corepresentation in the group $\Gamma_{\q}$ label different phases, under the assumption of enough strongly projective corepresentations. 
\ermrk

%\brmrk
%However, we have one drawback in our proposed scheme: we can not yet justify strongly projective corepresentation by any physical argument. This will be an interesting open question for the physicists.
%\ermrk

 %Hence, let us  consider $SO_{q}(3)$- symmetry,  In the notation used before, suppose (without loss of generality) that $N^{-1}_{[k]} =U_{\frac{m}{2}} (1 \ot q_1),$ $M_{[k]} =U_{\frac{n}{2}} (1 \ot q_2)$ for some $q_1,q_2 \in SO_q(3)$. It follows that
 %\blmma
 %\label{phases}
 %Either both $m,n$ are even or both are odd.
 %\elmma 
% \begin{proof} Suppose $m$ is odd and $n$ is even. Then from the fusion rules of $SU_q(2)$ corepresentations, we see that the matrix coefficients of $N^{-1}_{[k]} \ot M_{[k]}$ are linear combinations of $U_{\frac{l}{2}} $ for odd $l$'s only. But on the other hand, $N^{-1}_{[k]} \ot M_{[k]}$ is also equal to 
 %$W (1 \ot q)$ for unitary corepresentation $W$ and a unitary element $q$ of $SO_q(3)$,  hence its matrix coefficients must be linear combinations of $U_{\frac{p}{2}} for only even $p$. As the matrix coefficients  of different irreducibles are linearly independent, this gives a contradiction.
% \end{proof} 
% Recall that $U_{\frac{m}{2}}$ for odd $m$, i.e. half-integers, correspond to the nontrivial element in 
%$\Gamma_{SO_q(3)} \cong Z_2,$  and for even $m$'s (genuine integers) they correspond to the trivial element. Thus, the statement of Lemma \ref{phases} implies that either all the  there will be exactly 2 distinct topological phases for $SO_{q}(3)$- symmetric $1-D$ spin systems. 
 
 To give a concrete $SO_{q}(3)$-symmetric model, we consider the set up of \cite{Phy2} and in particular example 6.3. In this case, if we restricted ourselves to group symmetry only, we would have symmetry given by  only  the maximal classical  subgroup of $SO_{q}(3),$  i.e. a maximal quantum subgroup  isomorphic to $C(G)$ for some group $G$. However, it is known that 
  for all $q \ne 1$,  $G \cong S^1$ (see, for example, Proposition 3.13 of \cite{kacthesis}. As $S^1$ does not have any nontrivial $2$-cocycle,   we conclude that  for $q \neq 1,$ consideration of classical group symmetry would not reveal the non-trivial topological phase which is there when we consider the full quantum symmetry of the model. This clearly illustrates the potential advantage of allowing quantum group symmetry.
%It proves that $\Gamma$ is different from $H^{2}_{uinv}(\q, S^1).$ One can also observe that there is an injective map from $\Gamma$ to $H^{2}_{uinv}(\q, S^1).$ 

As a particular example of a one dimensional model with a global $SO(3)$ symmetry and a known symmetry protected topological ground state, we consider the Affleck-Kennedy-Lieb-Tasaki (AKLT) model \cite{aklt_original}. The model consists of quantum spins $S = 1$ arranged in a one-dimensional chain, with a Hamiltonian that reads
\begin{equation}
H = \sum_{i=1}^{N}\frac{1}{2}\vec{S}_{i}\cdot\vec{S}_{i+1} + \frac{1}{6}\sum_{i=1}^{N}\left(\vec{S}_{i}\cdot\vec{S}_{i+1}\right)^{2}
\label{aklt_hamilt}
\end{equation}

where the $\vec{S}_{i}$ are the spin operators on the site $i$ on the chain. The AKLT model is known to have a unique gapped ground state under periodic boundary conditions, while it harbors four possible edge states (an edge group $Z_{2}\times Z_{2}$ under open boundary conditions, which are topological in nature. The edge states are protected by the $SU(2)$ (or, $SO(3)$) symmetry of the Hamiltonian, and can be visualized as free spin-$\frac{1}{2}$ degrees of freedom at the two edges of the chain. Thus the AKLT model is a prototypical example of a symmetry protected topological state (SPT)\cite{zeng_chen_zhou_wen}. 

The $q$-deformation of the AKLT chain follows by rewriting the generators of the $SU(2)$ group, the spin operators in the spin-$1$ representation $\hat{S_{x}}$, $\hat{S_{y}}$ and $\hat{S_{z}}$ after $q$-deformation. The methods of doing this are quite standard, and here, we follow the discussion in \cite{simon}. Under the deformation, the AKLT Hamiltonian takes the more complicated form\cite{batchelor1990,franke_quella_2024}
\begin{eqnarray}
H_{q}/b & = & \sum_{i=1}^{N-1}c\vec{S}_{i}\cdot\vec{S}_{i+1} + \vec{S}_{i}\cdot\vec{S}_{i+1} + \frac{1}{2}\left(1-c\right)\left(q + q^{-1} - 2\right)\hat{S}_{i}^{z}\hat{S}_{i+1}^{z}\nonumber\\
&  & +\frac{1}{4}\left(1+c\right)\left(q - q^{-1}\right)\left(\hat{S}_{i+1}^{z} - \hat{S}_{i}^{z}\right)^{2}+ \frac{1}{4}c\left(1-c\right)\left(q + q^{-1} - 2\right)^{2}\left(\hat{S}_{i}^{z}\hat{S}_{i+1}^{z}\right)^{2}\nonumber\\
&  & +\frac{1}{4}c\left(1+c\right)\left(q + q^{-1}\right)\left(q + q^{-1} - 2\right)\hat{S}_{i}^{z}\hat{S}_{i+1}^{z}\left(\hat{S}_{i+1}^{z} - \hat{S}_{i}^{z}\right)\nonumber\\
&  &  +\frac{1}{4}\left(c-3\right)\left[c-1 +\frac{1}{2}\left(c+1\right)^{2}\hat{S}_{i}^{z}\hat{S}_{i+1}^{z} + 2\left(c-\frac{1}{8}\left(c+1\right)^{2}\right)\left(\left(\hat{S}_{i+1}^{z}\right)^{2} + \left(\hat{S}_{i}^{z}\right)^{2}\right)\right]\nonumber\\
&  & + \left(c-1\right) + \frac{1}{2}c\left(q^{2} - q^{-2}\right)\left(\hat{S}_{i+1}^{z} - \hat{S}_{i}^{z}\right)
%\label{qdeformedaklt}
\end{eqnarray}
where $c = 1 + q^{2} + q^{-2}$ and $b = \left[c\left(1-c\right)\right]^{-1}$. It is clear that the Hamiltonian matches the usual AKLT form in the limit $q\rightarrow 1$. When $q \neq 1$, the classical $SO(3)$ symmetry of the original AKLT model is lost, as is  evident from the Hamiltonian of the $q$-deformed model. Since the topological nature of the gapped ground state of the AKLT model (as well as the existence of the boundary or edge states) is crucially dependent on the classical $SO(3)$ symmetry, we would conclude that the $q$-deformed AKLT model will have no symmetry protected topological ground state. However, considering our formalism, we will predict that the $q$-deformed AKLT model will have \emph{two} distinct topological phases, one a topologically trivial phase and a topologically non-trivial phase, since $\Gamma_{SO_{q}(3)} = Z_{2}$. This result thus provides the evidence that by considering the cohomology of the underlying quantum group of a model, we can uncover new topological phases that will be absent if we restrict our attention to the classical symmetry group only. The elucidation of the actual nature of the topological phase in the $q$-deformed AKLT model will form the basis of future work. 
\paragraph{Acknowledgement}
D. Goswami is partially supported by JC Bose National Fellowship given by SERB, Govt. of India.
K. Maity is grateful to the Junior Research Fellowship in Mathematics granted to him by the Indian Statistical Institute, Kolkata.

\end{document}